\documentclass[pdflatex,sn-mathphys-num]{sn-jnl}% Math and Physical Sciences Numbered Reference Style
\usepackage{
    amsfonts,
    amsmath,
    amssymb,
    amsthm,
    booktabs,
    csvsimple,
    dsfont,
    graphicx,
    hyperref,
    mathrsfs, 
    mathtools,
    manyfoot,
    pgfplots, 
    suffix,
    textcase,
    textcomp,
    tikz,
    xcolor,
    eqparbox
} 

\pgfplotsset{compat=1.18}

\usepackage[title]{appendix}%

\usepackage{placeins}

\usetikzlibrary{fadings} 
\usetikzlibrary{positioning,pgfplots.groupplots}

\DeclarePairedDelimiter\norm{\lVert}{\rVert}

\newcommand*{\Lpnorm}[3][2]{\norm{#3}_{L^{#1}(#2)}}
\WithSuffix\newcommand\Lpnorm*[3][2]{\norm*{#3}_{L^{#1}(#2)}}

\newcommand*{\Wnorm}[4][2]{\norm{#4}_{W^{#2,{#1}}(#3)}}
\WithSuffix\newcommand\Wnorm*[4][2]{\norm*{#4}_{W^{#2,{#1}}(#3)}}

\newcommand*{\lpnorm}[3][2]{\norm{#3}_{\ell^{#1}(#2)}}
\WithSuffix\newcommand\lpnorm*[3][2]{\norm*{#3}_{\ell^{#1}(#2)}}

\newcommand*{\Hnorm}[3][-1]{\norm{#3}_{H^{#1}(#2)}}
\WithSuffix\newcommand\Hnorm*[3][-1]{\norm*{#3}_{H^{#1}(#2)}}

\newcommand{\NN}{{\cal N}}
\newcommand{\patch}{{\cal P}}
\newcommand{\at}[2][]{#1|_{#2}}

\theoremstyle{thmstyleone}%
\newtheorem{theorem}{Theorem}%  meant for continuous numbers
\newtheorem{lemma}{Lemma}%  meant for continuous numbers
\theoremstyle{thmstyletwo}%
\newtheorem{remark}{Remark}%

\theoremstyle{thmstylethree}%
\newtheorem{definition}{Definition}%

\begin{document}

\title[
    Error analysis for hybrid finite element/neural network discretizations
]{
    Error analysis for hybrid finite element/neural network discretizations
}

%%=============================================================%%
%% GivenName	-> \fnm{Joergen W.}
%% Particle	-> \spfx{van der} -> surname prefix
%% FamilyName	-> \sur{Ploeg}
%% Suffix	-> \sfx{IV}
%% \author*[1,2]{\fnm{Joergen W.} \spfx{van der} \sur{Ploeg} 
%%  \sfx{IV}}\email{iauthor@gmail.com}
%%=============================================================%%

\author[1]{\fnm{Uladzislau} \sur{Kapustsin}}\email{uladzislau.kapustsin@ovgu.de}

\author[1]{\fnm{Utku} \sur{Kaya}}\email{utku.kaya@ovgu.de}
% \equalcont{These authors contributed equally to this work.}

\author[2]{\fnm{Johannes} \sur{Pfefferer}}\email{johannes.pfefferer@unibw.de}

\author[1]{\fnm{Thomas} \sur{Richter}}\email{thomas.richter@ovgu.de}

\affil[1]{\orgname{Otto-von-Guericke-Univerist\"at Magdeburg}, \orgaddress{\street{Universita\"tsplatz 2}, \city{Magdeburg}, \postcode{39106}, \country{Germany}}}

\affil[2]{\orgname{Universit\"at der Bundeswehr M\"unchen}, \orgaddress{\street{Werner-Heisenberg-Weg 39}, \city{Neubiberg}, \postcode{85579}, \country{Germany}}}

%%==================================%%
%% Sample for unstructured abstract %%
%%==================================%%

\abstract{
    We analyze a hybrid finite element/neural network method for predicting solutions of partial differential equations. The methodology is designed for updating a coarse finite element solution, which is obtained efficiently by a standard finite element solver, with fine scale fluctuations from neural networks. The neural network is applied in a local manner in such a way that the coarse solution is split into \emph{patches}, unions of some few elements, and the network works independently on each patch. This supports generalization by design as the network never sees the full solution or even the domain. 
    In previous studies with applications to 3D Navier-Stokes flows we have shown accelerations up to 30 times when compared to usual Newton-multigrid solvers. The focus of this study is on the theoretical analysis of the hybrid method for the Poisson problem. 
    Key observation is the interplay between the fineness of the coarse finite element mesh, the size of the patches, the fineness of the fine prediction mesh and the quality and quantity of the training data. We provide the a priori error analysis of the method together with the stability analysis of the neural network. The numerical experiments confirm the capability of the network predicting fine finite element solutions. We also illustrate the generalization of the method to problems where test and training domains differ from each other.
}

\keywords{Neural networks, PDE approximation, Finite element method}
\pacs[MSC Classification]{65Y10, 65N12, 65N30}

\maketitle

\section{Introduction and motivation}

The use of neural networks to approximate solutions to partial differential equations has made considerable progress in recent years. In particular the class of Physics Informed Neural Networks (PINN)~\cite{raissi2019}, like the Deep Ritz method, has been investigated intensely~\cite{weinan2018}. Especially for high dimensional problems~\cite{pmlr-v134-lu21a,blechschmidt2021ways,schwab2023} or parameter dependent partial differential equations~\cite{belbuteperes2021} these approaches offer a new point of view and promise a substantial increase in efficiency. There are further methods in the literature which exploit the weak formulations of the PDEs  \cite{kharazmi2019} or incorporate finite element (FE) spaces in the loss function \cite{badia2023,Meethal2023,MITUSCH2021110651}. See~\cite{Tanyu2023} for a review.

However, when standard problems such as three dimensional fluid mechanics are considered, neural network based approaches have to compete with highly sophisticated and refined classical discretization methods. Finally, the training of the network often remains the crucial problem. Even though the theoretical approximation properties of neural networks may be superior, the efficient solution of the associated optimization problems is still an open problem. On the other hand, there are highly efficient Newton-Krylow space methods, possibly with optimal multigrid methods for preconditioning. 

The main drawback of simple PINN's is the need to re-train when the problem parameters change. DeepONet~\cite{lu2021} is one approach to this issue. Here, instead of training one net to represent a specific problem, two nets are used to represent possible solutions as well as the solution operator. But even this approach could not yet prove the efficiency and accuracy of established methods.

An alternative approach we follow is to combine established methods for coarse representation of a solution with a neural network to resolve fine scales whose representation is often not possible.
This approach has many potential applications, for example in fluid mechanics, where finite element or finite volume methods can reproduce the coarse structure with great accuracy while respecting conservation principles, but simultaneous resolution of fine-scale turbulent processes is often not possible or would simply be too expensive.
Corresponding approaches allow to reduce modeling as well as compute efforts significantly, while providing a high robustness~\cite{A41}. By combining data with classical numerical simulations, purely neural network driven methods can sometimes be outperformed~\cite{A42}.
In the context of PDE solving there are examples of machine learning architectures that take inspirations from numerical and mathematical ideas~\cite{A32,A55,A56}. These approaches mostly mimic discretisation concepts combined with iterative solving strategies in the underlying machine learning architectures. 
 A super-resolution methodology in this direction was introduced in  \cite{jiang2020}. With the \emph{Deep Neural Network Multigrid Solver} (DNN-MG)~\cite{Margenberg2021,Margenberg2022} we have introduced a concept which embeds a neural network fluently into a multigrid hierarchy, solves the coarse grid levels directly, e.g. with a finite element multigrid method, and predicts the corrections on fine grid levels locally by a neural network. Particularly in~\cite{Margenberg2023} we have demonstrated a 30 times acceleration in applications to nonstationary 3D Navier-Stokes simulations, when compared to state of the art parallel Newton-multigrid solvers. This is mostly due to the highly efficient evaluation of neural networks on accelerator cards close to peak performance. In a 3d Navier-Stokes benchmark presented in~\cite{Margenberg2023} we considered networks with up to $5\,000\,000$ parameters to enrich a finite element solution with just $400\,000$ unknowns and still found that the time for network evaluation was less than $3\%$ of the overall runtime. Modern accelerator hardware is tuned for neural network evaluation such the number of network parameters can be orders of magnitude higher than the number of finite element unknowns without introducing substantial overhead.

So far, these efforts have taken a computational approach and are mostly strictly heuristic and driven by the usual optimisations on the neural network side. A very first analysis has been done in our preliminary work~\cite{Kapustsin2023}. Besides the restriction to the stationary Laplace equations, however, we have not considered the decisive feature of DNN-MG: Instead of a local approach that applies corrections in terms of multiscale methods to patches from a few cells, we investigated the global application of the network. However, this approach is completely at odds with the desired generalisation through domain independence. 

Here, the neural network intervenes only locally: the grid is decomposed into patches, e.g., in 2d a range of $p\times p$ ($p\in\mathbb{N}$ is small here, usually less than four) elements, and on each of these patches the correction to a finer solution is pulled from the network. This approach was shown to increase efficiency for standard flow problems compared to established methods. Furthermore, the local design, i.e. the application of a mesh to all patches, allows a very good generalizability. The DNN-MG solver can be regarded as a numerical solution method in the sense of a domain decomposition method rather than an approximation method. The network never sees the whole solution but always only small sections. This local approach has similarities to non-overlapping domain decomposition methods~\cite{Chan_Mathew_1994}, however, with subdomains of size $H\to 0$ when $h\to 0$. There are also similarities to multiscale methods such as the heterogeneous multiscale method, where local problems are solved to resolve micro features~\cite{abdulle}.

The mathematical analysis of PINNs is already well advanced. In particular, the aspect of approximation with neural networks is very well established, starting with the universal approximation theorem in the late 80s~\cite{barron1993,cybenko1989}. Especially for the application to partial differential equations relevant results are available~\cite{Guehring2020,pmlr-v134-lu21a,muller2022} and optimal estimates in $W^{n,p}$-spaces are known, see~\cite{DeVore2021} for a review of neural network's approximation properties. Of the PINNs, the Deep Ritz method~\cite{weinan2018} in particular is well studied. It is based on the direct approximation of the energy functional (in the case of symmetric differential operators, e.g., Laplace or Stokes) with neural networks and Monte Carlo integration. Here, quite comprehensive a priori~\cite{muller2022} as well as a posteriori error estimation are available~\cite{minakowski2023}.

The goal of this work is to theoretically investigate hybrid approaches that enrich a finite element solution on coarse grids with fine scale fluctuations from a neural network in terms of the DNN-MG method. We restrict ourselves to the simple linear Poisson equation and give the complete a priori error analysis of the hybrid method. The efficiency in relevant applications has already been documented in \cite{Margenberg2023} using the example of the 3D Navier-Stokes equations. The aim here is to take a first decisive step towards the rigorous analysis of such hybrid methods.

In the next section we will briefly describe the finite element discretization and introduce some notation. Section~\ref{sec:hybrid} then introduces the hybrid approximation method and describes the training of the neural networks. We start the analysis of the method with the version where global network updates are used and then refine the a priori analysis for the local neural network updates. Numerical demonstrations follow in Section~\ref{sec:num}.

\section{Preliminaries}
We start with describing the model problem and then introduce the notation and finite element discretization. Let $\Omega\subset\mathbb{R}^d, \; d \in \{2,3\}$ be a domain with polygonal boundary. For $f\in H^{-1}(\Omega)$ let $u\in H^1_0(\Omega)$ be the weak solution to the Poisson equation
\begin{equation}
  -\Delta u = f, \quad 
  u \vert_{\partial \Omega} = 0.
  \label{eq:poisson equation}
\end{equation}
\subsection{Finite element discretization}\label{sec:fem}
By $\Omega_h$ we denote a finite element mesh of the domain $\Omega$, where $h$ stands for the diameter of the largest cell. For simplicity we assume that $\Omega_h$ is a quadrilateral or triangular mesh in 2D or a hexahedral or tetrahedral mesh in 3D, respectively, that satisfies the common assumptions of structural regularity and shape regularity~\cite{Ciarlet}. We assume that there is a hierarchy of finite element meshes 
\begin{equation}\label{triahierarchy}
  \Omega_{h_{\mathcal{P}}} 
  := 
  \Omega_{h_0} 
  \preccurlyeq \Omega_{h_1}
  \preccurlyeq \cdots
  \preccurlyeq \Omega_H:=\Omega_{h_{L'}}
  \preccurlyeq \cdots
  \preccurlyeq \Omega_{h_L} 
  =: 
  \Omega_h,
\end{equation}
where we denote by $\Omega_{h_{l-1}} \preccurlyeq \Omega_{h_l}$ that each element of the fine mesh $T \in \Omega_{h_l}$ originates from the uniform refinement of a coarse element $T'\in \Omega_{h_{l-1}}$, for instance, uniform splitting of a quadrilateral or triangular element into four and of a hexahedral or tetrahedral element into eight smaller ones, respectively. We further assume that the meshes $\Omega_H$ and $\Omega_h$ are quasi-uniform. This assumption will simplify the theory. However, we already point out that quasi-uniformity of $\Omega_H$ and $\Omega_h$ on each element $E\in \Omega_{h_\patch}$ might be sufficient.

On $\Omega_h$ let $V_h$ be the space of basis functions of degree $r \ge 1$ satisfying the homogeneous Dirichlet condition on the boundary $\partial\Omega$
\begin{equation*}
  V_h
  :=
  \left\{
  \phi \in C(\bar \Omega), \; 
  \phi \at{T} \in P^{(r)}(T) \; \forall T \in \Omega_h, \;
  \phi \at{\partial \Omega} = 0
  \right\}
\end{equation*}
where, on triangular meshes, $P^{(r)}(T)$ is the space of polynomials of maximum degree $r$ and, on quadrilateral meshes, $P^{(r)}(T)$ is the tensor-product space of polynomials of maximum degree $r$ in every coordinate direction. These spaces are nested, i.e.,
\[
V_{h_{l-1}} \subset V_{h_l}, \quad l = 1, \dots, L.
\]
Then, $u_h \in V_h$ is the finite element solution to
\begin{equation}\label{fem}
  (\nabla u_h, \nabla\phi_h) 
  = 
  \langle f, \phi_h\rangle \quad \forall \phi_h \in V_h,
\end{equation}
with $\langle\cdot,\cdot\rangle$ denoting the duality pairing between $H^{-1}(\Omega)$ and $H^1_0(\Omega)$.
For two right hand sides $f, g \in H^{-1}(\Omega)$ and the corresponding finite element solutions $u_h^f, u_h^g \in V_h$ it holds
\begin{equation}\label{femstability}
  \Lpnorm{\Omega}{\nabla (u_h^f - u_h^g)} 
  = 
  \Hnorm{\Omega}{f-g}.
\end{equation}

Given sufficient regularity of the right hand side, namely $f \in H^{r-1}(\Omega)$ (using for simplicity the notation $H^0(\Omega) := L^2(\Omega)$), we forthcoming assume that $u$ belongs to $H^{r+1}(\Omega)$. In addition, we assume convexity of the domain $\Omega$. As a consequence, the standard a priori estimate
\begin{equation}\label{femerror}
  \Lpnorm{\Omega}{u-u_h}
  + 
  h \Lpnorm{\Omega}{\nabla (u-u_h)}
  \leq 
  c h^{r+1} \Hnorm[r-1]{\Omega}{f} 
\end{equation}
holds. For instance, the regularity $u\in H^{r+1}(\Omega)$ is granted for $r=1$ if $\Omega$ is convex. For smaller interior angles, we can even guarantee higher regularity, see \cite{Grisvard1985}.
\subsection{Notation}
Over subdomains $\omega \subseteq \Omega$ the $L^2-$ and $H^s-$norms are denoted by $\|\cdot\|_{\omega}$ and $\|\cdot\|_{s,\omega}$, respectively. The corresponding inner products are denoted by $(\cdot,\cdot)_{\omega}$ and $(\cdot,\cdot)_{s,\omega}$, respectively. We omit the index $\omega$, if the norm or the inner product is considered on the whole domain $\Omega$. 
With $\|\cdot\|_2$ we denote the Euclidean norm for vectors and the spectral norm for matrices. With $X^h_{\omega}$ and $X^H_{\omega}$ we denote the nodes of the meshes $\Omega_h$ and $\Omega_H$ that lie in the subdomain $\overline{\omega}$, respectively. Moreover, for $v \in C(\overline{\omega})$ we define 
\begin{equation}\label{eq:l2_norm}
  \| v \|_{\ell^2(\omega)} :=\Big( \sum \limits_{x \in X_{\omega}^h} v(x)^2 \Big)^{\frac1{2}}. 
\end{equation}

\begin{definition}[Patch]\label{def:patch}
  A patch $\patch$ is defined to be a subdomain that is geometrically identical to the union $M$ of certain cells of the coarse mesh $\Omega_{h_{\mathcal{P}}}$. For simplicity, we assume that for each patch there holds $|\patch|=\mathcal{O}(h_\patch)$.
  We also exploit that a patch is not only identified with the degrees of freedoms of this element but also by the cells assembling it, $\patch = \{T\in\Omega_h \;: T \subset M  \}$. The set of all patches will be denoted by $U_\patch$.
\end{definition}
By $V_h(\patch)$ we denote the local finite element subspace
\[
V_h(\patch) := \operatorname{span} \left\{
\phi_h \at{\patch}, \; \phi_h \in V_h
\right\}. 
\]
By $R_\patch : V_h \to V_h(\patch)$ we denote the restriction to the local patch space, defined via
\[
R_\patch(u_h)(x) = u_h(x) \quad \forall x\in \patch.
\]
By $P_\patch : V_h(\patch) \to V_h$ we denote the prolongation
defined by
\[
P_\patch(u_\patch)(x) = \begin{cases}
  \frac{1}{n(x)} u_\patch(x) & x \in \patch, \\
  0  & \text{otherwise},
\end{cases}
\]
where $n(x) \in \mathbb{N}$ is the number of patches that contain the degree of freedom $x$.

\section{Hybrid finite element neural network discretization}\label{sec:hybrid}
Consider two finite dimensional spaces $V_H$ and $V_h$ that are built on coarse and fine meshes $\Omega_H$ and $\Omega_h$, respectively.
The idea of the paper is to determine an approximate solution on the coarse mesh $\Omega_H$ with the finite element method and then to obtain the fine mesh fluctuations in forms of neural network updates.
In other words, we seek hybrid solutions $u_{\NN} \in V_h$ which are found by augmenting $u_H \in V_H$ with a neural network update $w_\NN \in V_h$, i.e.
\[
u_\NN := u_H + w_\NN. 
\]

\begin{figure}[t]
  \begin{center}
    \includegraphics[width=0.5\textwidth]{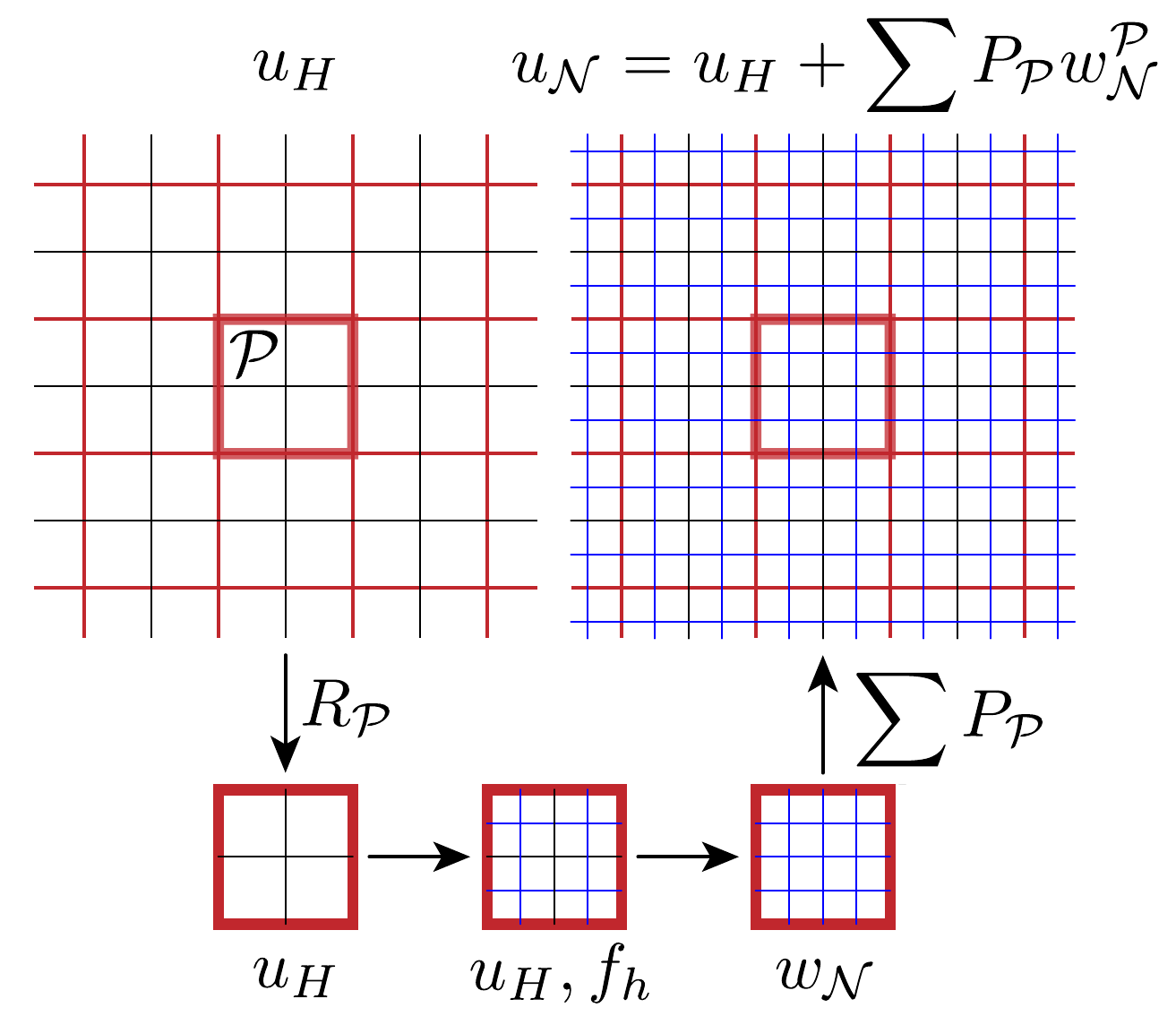}
    \caption{Illustration of the hybrid solver. The finite element solution $u_H$ is approximated on the coarse mesh $\Omega_H$ (in black).
    A patch mesh $\Omega_{\cal P}$ which has the same level or is even coarser than $\Omega_H$ (in red, here, one even coarser) combines elements from the coarse mesh. 
    On each patch the solution is locally extracted. Together with the fine mesh right hand side information on $\Omega_h$ (shown in blue) it is the input of a neural network. The output $w_{\cal N}$ is the local correction towards an improved solution and is gathered on the fine global mesh $\Omega_h$. }
    \label{fig:method_illustration}
  \end{center}
\end{figure}

The neural network predicts the finite element coefficients of the correction $w_\NN$ and as input it receives all data that is available: the coarse solution $u_H \in V_H$ and the problem data, i.e. the source term $f$.
Moreover, these updates are to be obtained locally in such a way that the
network is not acting on the aforementioned data on the whole domain $\Omega$. 
Instead the network acts on the patches $\patch$ separately in order to obtain the values of an update $w_\NN \at{\patch}$ by providing the coefficient vector for its representation in $V_h(\patch)$.
An illustration of the local updates is given in Fig.~\ref{fig:method_illustration}. 

\begin{definition}[Hybrid solution] \label{def:hybrid solution}
  Given a subdivision of the domain $\Omega$ into a set of patches, a neural network function $\NN(\cdot)$ acting on each patch, the coarse mesh solution $u_H\in V_H$ and the right hand side $f$, we define the hybrid solution $u_\NN \in V_h$ as 
  \begin{equation}\label{eq:hybridsolution}
  u_\NN  
  := 
  u_H 
  + 
  w_\NN
  =
  u_H 
  +
  \sum\limits_{\mathclap{\patch\in U_\patch}} 
  P_{\patch} w_\NN^{\patch},
  \end{equation}
  where 
  $w_\NN^{\patch} =\sum_{i=1}^{N_{dof}} \NN(\mathrm{x}_\patch)_i  \phi^\patch_i$. Here, $\NN(\mathrm{x}_\patch)_i$ is the $i$-th output of the network, $\phi_i^\patch ,\; i \in \{ 1,\cdots,N_{dof}\}$, are the basis functions of $V_h(\patch)$ and $\mathrm{x}_\patch$ are network inputs corresponding to a patch $\patch$ defined as the vector $\mathrm{x}_\patch=(\mathrm{x}_\patch^{u_H},\mathrm{x}_\patch^f)$, where
\begin{equation} \label{eq:network inputs}
\mathrm{x}_{\patch}^{u_H} = \begin{pmatrix}  u_H(x) \end{pmatrix}_{x \in X^H_\patch},
\quad \mathrm{x}_{\patch}^{f}  = \begin{pmatrix} f(x) \end{pmatrix}_{x \in X^h_\patch},
\end{equation}
and the order of degrees of freedom is consistent across all patches and elements.
\end{definition}
As networks we will only consider fully connected multilayer perceptrons:

\tikzset{%
  every neuron/.style={
    circle,
    draw,
    minimum size=.5cm
  },
  neuron missing/.style={
    draw=none, 
    scale=2,
    text height=0.333cm,
    execute at begin node=\color{black}$\vdots$
  }
}

\begin{figure}[t]
  \centering
\begin{tikzpicture}[
    neuron/.style={circle, draw, minimum size=8mm},
    input/.style={neuron, fill=blue!30},
    hidden/.style={neuron, fill=orange!40},
    output/.style={neuron, fill=green!40},
    conn/.style={->, thick}
]

% =====================
% Input layer (partial)
% =====================
\node[input] (I1) at (0,-1) {};
\node[input] (I2) at (0,-2) {};
\node at (0,-3) {$\vdots$};
\node[input] (I4) at (0,-4) {};

\node[above=6mm of I1] {Input Layer};

% Input data brace
\draw[decorate, decoration={brace, amplitude=6pt}]
  (-0.7,-4.4) -- (-0.7,-0.6)
  node[midway, xshift=-12pt, rotate=90] {Input Data};

% =====================
% Hidden layer 1 (partial)
% =====================
\node[hidden] (H11) at (2,-0.5) {};
\node[hidden] (H12) at (2,-1.5) {};
\node at (2,-2.8) {$\vdots$};
\node[hidden] (H15) at (2,-4.5) {};

% =====================
% Hidden layer 2 (partial)
% =====================
\node[hidden] (H21) at (4,-0.5) {};
\node[hidden] (H22) at (4,-1.5) {};
\node at (4,-3.2) {$\vdots$};
\node[hidden] (H24) at (4,-4.5) {};

% Hidden layers brace
\draw[decorate, decoration={brace, mirror, amplitude=6pt}]
  (1.3,-5.0) -- (4.7,-5.0)
  node[midway, yshift=-12pt] {Hidden Layers};

% =====================
% Output layer (partial)
% =====================
\node[output] (O1) at (6,-1) {};
\node at (6,-3) {$\vdots$};
\node[output] (O2) at (6,-2) {};
\node[output] (O4) at (6,-4) {};

\node[above=6mm of O1] {Output Layer};

% Output data brace
\draw[decorate, decoration={brace, mirror, amplitude=6pt}]
  (6.7,-4.4) -- (6.7,-0.6)
  node[midway, xshift=12pt, rotate=270] {Output Data};

% =====================
% Connections (representative)
% =====================
\foreach \i in {1,2,4}
  \foreach \j in {11,12,15}
    \draw[conn] (I\i) -- (H\j);

\foreach \i in {11,12,15}
  \foreach \j in {21,22,24}
    \draw[conn] (H\i) -- (H\j);

\foreach \i in {21,22,24}
  \foreach \j in {1,2,4}
    \draw[conn] (H\i) -- (O\j);

\end{tikzpicture}
  \caption{Multilayer perceptron}
\end{figure}

\begin{definition}[Multilayer perceptron]\label{def:MLP}
  A function ${\cal N} : \mathbb{R}^{N_{0}} \to \mathbb{R}^{N_{L}}$ of the following form
  \begin{equation}
    \NN := l_{L} \circ \sigma \circ l_{L-1} \circ\dots \circ \sigma \circ l_1
    \label{eq:mlp}
  \end{equation}
  where 
  \begin{equation}\label{def:layer}
    l_i(x) =  W_i x + b_i, \; i = 1, \dots ,{L} ,\; W_i \in \mathbb{R}^{N_{i-1} \times N_{i}},\; b_i \in \mathbb{R}^{N_i}
  \end{equation}
  and $\sigma : \mathbb{R} \to \mathbb{R}$ is called a \emph{multilayer perceptron (MLP)} of \emph{depth} $L \in \mathbb{N}$. 
  Here $\sigma$ is called \emph{activation function}, $l_i$ - \emph{layers}, $W_i$ - \emph{weights} and $b_i$ - \emph{biases}.
\end{definition}

We have studied the effect of the neural network architecture on the numerical performance of the hybrid solver in previous work~\cite{Margenberg2021,margenberg2021deepneuralnetworksgeometric,DNNnetworks2025} and found that the significance of the proper design is more important for stability and robustness than for accuracy.

\subsection{Training of the neural network}

Now, let us explain how we pre-train a neural network $\NN$.
First of all, it is necessary to select a set of training problems. We generate the training data by selecting a  set ${\cal{F}}\subset H^{-1}(\Omega)$ of right hand side functions.
Then, we solve the Poisson equation for each $f_i \in \cal{F}$ for both coarse and fine meshes $\Omega_H$ and $\Omega_h$, respectively. The solutions are denoted by $u_H^i$ and $u_h^i$, respectively.
The vector of input data \[\label{eq:x_patch^i}\mathrm{x}_{\patch}^i=\left(\mathrm{x}_\patch^{u_H^i},\mathrm{x}_\patch^{f_i}\right)\]
for the training of the neural network is provided for each patch $\patch\in U_\patch$ as described in \eqref{eq:network inputs}.
The corresponding output data $\mathrm{y}^i_\patch\in \mathbb{R}^{N_L}$ is given as a vector defined as
\[
\mathrm{y}^i_{\patch} = \begin{pmatrix}
  u^i_h(x) - u^i_H(x)
\end{pmatrix}^T_{x \in X^h_{\patch}},
\]
which is simply a difference between fine and coarse solution on each fine mesh node $x$ belonging to the patch $\overline{\patch}$ and again the order of degrees of freedom is consistent across all patches and elements.
We would like to underline once more that the network update of the correction $w_\NN$ is given by the value of the network applied not globally to the whole domain $\Omega$, but locally to each patch $\patch$ (see Definition~\ref{def:patch}).

Once the training data is available, we need to solve the following optimization problem
\begin{equation}\label{minimize}
  \min\limits_{\substack{W_j\in\mathbb{R}^{N_{j-1},N_i},b_i\in\mathbb{R}^{N_j},\\ j \in \{0,\dots,L\}}} 
  \quad 
  \frac{1}{N_T N_P}
  \sum\limits_{\substack{i \in \overline{1, N_T} \\ \patch \in U_\patch}} \norm{\mathrm{y}^i_\patch -  \NN(\mathrm{x}_\patch^i)}^2_2   
\end{equation}
where $N_T$ is the number of source terms in the training set ${\cal F}$ and $N_P$ is the number of patches.
In order to solve this problem we use one of the stochastic gradient descent based methods. After the network has been trained, we can finally apply it to other problems that were not in the training data. For this, firstly we again construct input data in the same way as described above and compute network predictions. Then, we construct a complete network solution $u_\NN$ by summing up (and averaging) these predictions as described in Definition~\ref{def:hybrid solution}.

Considering everything we have developed so far from the finite element perspective, we note that each coarse solution $u_H \in V_H$ also belongs to $V_h$ and there it takes the form
\[
u_H = \sum\limits_{i=1}^{N_{dof}} U^i_h \phi^i_h,
\]
where  $\{\phi^i_h\}_{i=1}^{N_{dof}}$ is the basis of the fine finite element space $V_h$ and where $U^i_h\in\mathbb{R}$ are corresponding coefficients. As a consequence of the fact, that we update the coarse solution $u_H$ only on fine mesh nodes, see \eqref{eq:hybridsolution}, we can consider this whole procedure as a simple update of fine mesh coefficients $U^i_h$, i.e.
\[
u_\NN = \sum\limits_{i=1}^{N_{dof}}	(U^i_h + W^i_\NN) \phi^i_h \in V_h.
\]

\subsection{Stability of the neural network}\label{sec:stability}

We continue with a stability result, which will be used to further discuss the error estimates in Theorems \ref{thrm1} and \ref{thrm2} later on, i.e., we analyze the patch-wise error between $w^\patch_\NN$ and $w_\NN^{\patch, i}$, where the first term is defined by $$w_\NN^{\patch} =\sum_{j=1}^{N_{dof}} \NN(\mathrm{x}_\patch)_j  \phi^\patch_j,$$ see Definition~\ref{def:hybrid solution}, and the second one is analogously defined by $$w_\NN^{\patch,i} =\sum_{j=1}^{N_{dof}} \NN(\mathrm{x}_\patch^i)_j \phi^\patch_j,$$ see \eqref{eq:x_patch^i} for the definition of $\mathrm{x}_\patch^i$.  However, let us first collect a couple of observations:
There exist constants $c_1, c_2 >0$ independent of $h_\patch$ and $n_\patch$, which is the number of degrees of freedom in each patch, such that a simple scaling argument gives 
\[
c_1 \|v\|_\patch^2 \le \frac{h_\patch^d}{n_\patch} \|v\|_{l^2(\patch)}^2 \le c_2 \|v\|_\patch^2\quad\forall v\in V_h(\patch).
\]
The number of fine mesh coefficients in each patch $\patch$ scales like $n_\patch = {\cal O}(h^{-d} h_\patch^d)$, hence
\begin{equation}\label{eq:norms}
c_1 \|v\|_\patch^2 \le h^d \|v\|_{l^2(\patch)}^2 \le c_2 \|v\|_\patch^2\quad\forall v\in V_h(\patch).
\end{equation}
Together with the inverse inequality and \eqref{eq:norms}, for each $v \in V_h(\patch)$ it holds
\begin{equation}\label{eq:norms2}
\|\nabla v\|_{\patch}^2 \leq  ch^{-2} \|v\|_{\patch}^2
\leq c h^{d-2} \lpnorm{\patch}{v}^2.
\end{equation}

\begin{lemma}[Network stability]\label{lemma:network_stability}
Let $\NN$ be a $L(\NN)$-Lipschitz continuous neural network, i.e., the neural network Lipschitz continuously maps the input data to the output data with Lipschitz constant $L(\NN)$. Then it holds that
\begin{equation}
    \norm{\nabla (w^\patch_\NN  - w_\NN^{\patch,i})}_\patch 
    \leq c L(\NN) h^{-1} \left(
      \norm{u_H - u_H^{i}}_{\patch} 
      +  \norm{I_h(f - f_i)}_{\patch}
    \right),
\end{equation}
where $I_h (f-f_i)\in V_h(\patch)$ denotes the Lagrange interpolant of $f-f_i\in C(\bar\Omega)$ on the patch~$\patch$.
\end{lemma}

\begin{proof}
We start by applying inequality \eqref{eq:norms2} to $\norm{\nabla (w^\patch_\NN - w_\NN^{\patch, i})}_\patch$ and use the fact that the nodal values of network updates are equal to the values of network outputs, i.e.,
  \begin{equation}  
\begin{split}
    \norm{\nabla (w^\patch_\NN - w_\NN^{\patch, i})}_\patch
    & \le c  h^{d/2-1} \norm{w^\patch_\NN - w_\NN^{\patch, i}}_{\ell^2(\patch)} \\
    & = c  h^{d/2-1} \norm{\NN(\mathrm{x}_{\patch}) - \NN(\mathrm{x}^{i}_\patch)}_2.
    \end{split}
  \end{equation}
Then, by using the Lipschitz continuity of MLP we obtain
\begin{equation}\label{eq:network_stability_1}
  \begin{split}
    \norm{\nabla (w_\NN - w_\NN^{i})}_\patch 
    & \leq c  L(\NN)  h^{d/2-1} \norm{\mathrm{x}_\patch - \mathrm{x}^{i}_\patch}_2.
    % & \leq c \cdot L(\NN) \cdot h^{-1} \left(
    %   \norm{u_H - u_H^{i}}_{\patch} 
    %   +  \norm{f - f_i}_{\patch}
    % \right)
  \end{split}
\end{equation}
  By definition of $\mathrm{x}_\patch$ and $\mathrm{x}^{i}_\patch$ we have
  \begin{equation} 
    \norm{\mathrm{x}_\patch- \mathrm{x}^{i}_\patch}^2_2
    = 
    \lpnorm{\patch}{u_H - u^{i}_H}^2
    + 
    \lpnorm{\patch}{f - f_i}^2
  \end{equation}
and by applying inequality \eqref{eq:norms} we get
\begin{equation}\label{eq:network_stability_2}
  \norm{\mathrm{x}_\patch - \mathrm{x}^{i}_\patch}_2 
  \leq 
  c h^{-\frac{d}{2}} (\norm{u_H - u_H^{i}}_{\patch}
  +
  \norm{I_h(f- f_i)}_{\patch}
  ).
\end{equation}
By combining inequalities \eqref{eq:network_stability_1} and \eqref{eq:network_stability_2} we get the assertion.
\end{proof}

\begin{remark}[Lipschitz constant of MLP]
There exist various methods for estimating upper bounds of Lipschitz constants, e.g. using automatic differentiation \cite{NEURIPS2018_d54e99a6} or methods based on a relaxation to a polynomial optimization problem \cite{latorre2020}. The simplest bound \cite[Proposition 1]{NEURIPS2018_d54e99a6} can be obtained for an MLP with 1-Lipschitz activations in the following form 
\begin{equation}
    L(\NN) \leq \prod\limits_{i=1}^L \norm{W_i}_2.
\end{equation}
\end{remark}

\subsection{Error estimate using a global network update}
\label{sec:global_update}

We will start with a simplified setting, where there is only one patch that covers the complete domain, $U_\patch =\{\Omega\}$, which is possible as by construction there holds $\overline\Omega=\bigcup_{T\in\Omega_{h_\patch}}\overline{T}$, compare Definition~\ref{def:patch}. Here, the hybrid finite element solution is directly given by
\[
u_\NN =  u_H + w_\NN^{\Omega}, 
\]
see Definition~\ref{def:hybrid solution}. In order to indicate, that the Poisson equation is approximately solved with right hand side $f_i\in\mathcal{F}$ , we also use the notation $u_h^{i}$ and $u_H^{i}$ for the fine and coarse grid solution, respectively, and
\[
    u_\NN^{i} = u_H^{i} + w_\NN^{i}
\]
for the corresponding hybrid solution.

\begin{theorem}[A priori finite element error for the single-patch solution]\label{thrm1}
  Let $U_\patch=\{\Omega\}$ and let $\NN$ be the $L(\NN)$-Lipschitz continuous network trained on the training set ${\cal F}=\{f_1,\dots,f_{N_T}\}$ with $f_i \in C(\bar \Omega)$, $i=1,\dots,N_T$, such that the loss function \eqref{minimize} is reduced to the error order of $\epsilon^2$ for each training problem. For $f\in H^{r-1}(\Omega)\cap C(\bar \Omega)$, let $u\in H^{r+1}(\Omega)\cap H^1_0(\Omega)$ be the solution to the Poisson problem and $u_\NN\in V_h$ be the hybrid solution. It holds
\begin{align*}
  \|\nabla(u - u_\NN)\|
  &\le 
  c \Bigg(
    h^r \|f\|_{r-1}
    + 
    \min_{f_i \in {\cal F}} \Big\{
      \|f - f_i\|_{-1} \\
  &\qquad\qquad
      + L(\NN)  h^{-1} 
      \big(
        \|u_H - u_H^{i}\|
        + \|I_h(f - f_i)\|
      \big)
    \Big\}
    + h^{\frac{d}{2}-1} \epsilon
  \Bigg).
\end{align*}
\end{theorem}
\begin{proof}
  For arbitrary $f_i\in {\cal F}$ we split the error into
  \begin{multline*}\label{e1:1}
    \|\nabla (u-u_\NN) \| \le
    \underbrace{\|\nabla (u-u_{h}) \|}_{=(I)}+ 
    \underbrace{\|\nabla (u_{h}-u_{h}^{i}) \|}_{=(II)}+
    \underbrace{\|\nabla (u_{h}^{i}-u_\NN^{i}) \|}_{=(III)}+
    \underbrace{\|\nabla (u_\NN^{i}-u_\NN) \|}_{=(IV)}.
  \end{multline*}
  The fine mesh finite element error is estimated by 
  \begin{equation*}\label{e1:2}
    (I) = \|\nabla (u-u_h)\| \le c h^r \|f\|_{r-1}
  \end{equation*}
  according to \eqref{femerror}. Next, using \eqref{femstability}, the data approximation error is bounded by
  \begin{equation*}\label{e1:3}
    (II) = \|\nabla (u_h-u_h^{i})\| \le \|f-f_i\|_{-1}.
  \end{equation*}
  By \eqref{eq:norms2} and the assumption that the loss is reduced to the order of $\epsilon^2$ for each problem from the training set, we get a bound for the network approximation and optimization error
  \begin{equation*}\label{e1:4}
    (III) = \|\nabla (u_h^{i}-u_\NN^{i})\| 
    \le c h^{\frac{d}{2}-1} \|u_h^{i}-u_\NN^{i}\|_{\ell^2(\Omega)} \le c h^{\frac{d}{2}-1} \epsilon.
  \end{equation*}
  Finally, the generalization error of the network and a further error term depending on the richness of the data set remains:   
  \begin{equation*}\label{e1:5}
    \begin{aligned}
      (IV) = \|\nabla (u_\NN^{i}-u_\NN)\|
      &\le
      \|\nabla (u_H^{i}-u_H)\| + \|\nabla(w_\NN^{i} - w_\NN)\|\\
      &\le
      \|f-f_i\|_{-1} + \|\nabla(w_\NN^{i}-w_\NN)\|,
    \end{aligned}
  \end{equation*}
  see \eqref{femstability} for the last step.
  The last term can be bounded as follows using Lemma~\ref{lemma:network_stability}
\begin{equation}
\|\nabla(w_\NN^{i}-w_\NN)\| \leq cL(\NN) h^{-1} 
      \big(
        \|u_H - u_H^{i}\|
        + \|I_h(f - f_i)\|
      \big).
\end{equation}
 Combining the above gives the result.
\end{proof}

This lemma shows that the hybrid approach is able to reduce the error up to the accuracy of the fine finite element space $V_h$, the tolerance of the neural network training, the richness of the training data set ${\cal F}$ and the stability of the neural network that governs it ability to generalize beyond the training data set.

\subsection{Error estimates for local neural network updates}

We now tend to the discussion of a local application of the neural network, i.e., we assume that the set of patches $U_\patch$ coincides with the set of elements in $\Omega_{h_\patch}$, i.e. $U_\patch=\Omega_{h_\patch}$.
The global analysis in Section \ref{sec:global_update} requires $f\approx f_i$ for some $i\in\{1,\ldots,N_T\}$ in the whole domain $\Omega$ to have a good approximation quality of the hybrid solution. The localized approach will however allow us for more flexibility when it comes to choosing an approximation of the right hand side to one of the training data  $f\approx f_i$, as this can be done individually on each patch.

For the main result of this section we require some more notation: For each patch $\patch\in U_\patch$ we define an enlarged domain $\tilde\patch$ with $\patch\subset\subset\tilde\patch$ that could be the union of all patches $\patch'\in U_\patch$ that overlap with $\patch$, i.e.
\begin{equation}\label{enlarged}
  \tilde\patch := 
  \smashoperator[r]{\bigcup_{
      \substack{\patch'\in U_\patch\\
        \bar\patch'\cap\bar\patch\neq \emptyset}}}\,\patch'.
\end{equation}
For the distance between $\patch$ and $\partial\tilde\patch\setminus\partial\Omega$ it holds $d(\patch,\tilde\patch)=\operatorname{dist}(\patch,\partial\tilde\patch\setminus\partial\Omega)={\cal O}(h_\patch)$, where $h_{\patch}$ is the diameter of the patch $\patch$.
This distance is relevant for local error estimates.

Now we present the main result of this paper which we will prove later on. Note, we assume that all meshes are simplicial. This restriction comes from the fact that we use local finite element estimates from \cite{DemlowGuzmanSchatz2011}, which are only proven for simplicial meshes.

\begin{theorem}[A priori finite element error for the hybrid solution based on local patches] 
  \label{thrm2}
  Let $U_\patch=\Omega_\patch$ and let the meshes $\Omega_\patch$, $\Omega_H$ and $\Omega_h$ have simplicial structure.
  Let $\NN$ be the $L(\NN)$-Lipschitz continuous network trained on the training set ${\cal F}=\{f_1,\dots,f_{N_T}\}$ with $f_i \in H^{r-1}(\Omega)\cap C(\bar \Omega)$, $i=1,\dots,N_T$,
  such that the loss \eqref{minimize} is reduced below the error order of $\epsilon^2$ on each patch for each training problem. 
  For $f\in H^{r-1}(\Omega)\cap C(\bar \Omega)$, let $u\in H^{r+1}(\Omega)\cap H^1_0(\Omega)$ be the solution to the Poisson problem and $u_\NN\in V_h$ be the hybrid solution. It holds
  \begin{equation}\label{aprioriestimate}
    \begin{aligned}
    \norm{\nabla  (u - u_\NN)}
    &\leq c \Bigg[
        \left(h^r+H^{r+1} h_\patch^{-1}\right)  \norm{f}_{r-1}
        +\Big(\sum\limits_{\patch \in U_\patch} \min\limits_{i \in \overline{1,N_T}} \Big\{
            h_\patch^{-2} \Big( \norm{u_h^i - u_H^i}_{\tilde{\patch}}^2\\
            &\quad+ H^2 \norm{\nabla (u_h^i - u_H^i)}_{\tilde{\patch}}^2\Big)
            + \norm{\nabla(u_H - u_H^i)}_{\tilde{\patch}}^2 + \norm{f-f_i}_{-1, \patch}^2 \\
            & \quad + L(\NN)^2h^{-2} \Big(
                \norm{u_H - u_H^i}_\patch^2 + \norm{I_h(f-f_i)}_\patch^2
            \Big)
        \Big\} \Big)^{\frac{1}{2}} + N_\patch^\frac12 h^{\frac{d}{2}-1}\epsilon
    \Bigg].
    \end{aligned} 
  \end{equation}
\end{theorem}

The proof of this theorem is postponed to the end of this section. As preparation for the proof we introduce local auxiliary problems posed on the patches~$\tilde P$.
 
\begin{definition}[Local problems]\label{def:local}  
  Let $\patch\subset\subset \tilde\patch$ be a patch and $\tilde \patch$ a slightly enlarged domain with ${\cal O}(h_\patch)=d(\patch,\tilde\patch) \ge ch$ (for details concerning the constant $c$ we refer to \cite[Lemma 3.3]{DemlowGuzmanSchatz2011}), which is convex and matches the mesh $\Omega_{h_\patch}$. By $V_h(\tilde\patch)$, $V_H(\tilde\patch)$ and $V_{h,0}(\tilde\patch)$ we denote local finite element spaces, $V_{h,0}(\tilde \patch)$ having zero boundary data on $\partial\tilde\patch$. For $u_H\in V_H(\tilde\patch)$ we define $v_h\in V_{h,0}(\tilde\patch)$ via
  \begin{equation}\label{locproblem} 
    (\nabla (u_H+v_h),\nabla\phi_h)_{\tilde\patch} = (f,\phi_h)\quad\forall \phi_h\in V_{h,0}(\tilde\patch).
  \end{equation}
  Note, due to \eqref{triahierarchy} and \eqref{fem}, the solution $v_h$ also satisfies
  \begin{equation}\label{locproblem2}
    (\nabla v_h,\nabla\phi_H)_{\tilde\patch} = 0 \quad\forall \phi_H\in V_{H,0}(\tilde\patch).
  \end{equation}
\end{definition}
\begin{lemma}[Local problems]\label{lemma:local}
  Let $u_h\in V_h$ and $u_H\in V_H$ be solutions to \eqref{fem} with mesh parameters $h$ and $H$, respectively.
  For the solution to the local problem $v_h\in V_{h,0}(\tilde\patch)$ it holds
  \begin{equation}
    \|\nabla v_h \|_{\tilde\patch}+H^{-1}\|v_h\|_{\tilde\patch}\le c\|\nabla (u_h-u_H)\|_{\tilde\patch} \label{loc1est}.
  \end{equation}
\end{lemma}
\begin{proof}
  Testing~\eqref{locproblem} with $\phi_h=v_h$ and inserting $\pm (\nabla u_h,\nabla v_h)_{\tilde\patch}$ gives
  \[
  \|\nabla v_h\|^2_{\tilde\patch} =
  (f,v_h)_{\tilde\patch}    -(\nabla u_h,\nabla v_h)_{\tilde\patch}
  +(\nabla (u_h-u_H),\nabla v_h)_{\tilde\patch}.
  \]
  The first two terms on the right hand side sum up to zero using~\eqref{fem}, as $v_h$ can be extended to $\hat v_h\in V_h$ by zero outside of $\tilde\patch$. Estimating with Cauchy-Schwarz and dividing by $\|\nabla v_h\|_{\tilde \patch}$ gives the energy norm estimate in \eqref{loc1est}.
  Next, for $\psi\in L^2(\tilde\patch)$ we define the adjoint solution $z\in H^1_0(\tilde\patch)$
  \begin{equation}
    (\phi,\psi) = (\nabla \phi, \nabla z)\quad\forall\phi\in H^1_0(\tilde\patch).
  \end{equation}
  Due to the convexity assumption on $\tilde P$, we obtain $z\in H^2(\tilde\patch)$ by elliptic regularity. Moreover, a corresponding a priori estimates in combination with a standard scaling argument shows that
  \begin{equation}\label{eq:H2_aprior_patch}
    \|\nabla^2z\|_{\tilde\patch}\le c\|\psi\|_{\tilde\patch},
  \end{equation}
  where the constant $c>0$ is independent of $h_{\mathcal{P}}$.
  For $\phi:=v_h\in V_{h,0}(\tilde\patch)\subset H^1_0(\tilde\patch)$ it holds by means of~\eqref{locproblem2} for all $z_H\in V_{H,0}(\tilde \patch)$ that
  \[
  \begin{aligned}
    \left|(v_h,\psi)_{\tilde\patch}\right| &= \left|(\nabla v_h,\nabla z)_{\tilde\patch}\right|=\left|(\nabla v_h,\nabla (z-z_H))_{\tilde\patch}\right|\le \|\nabla v_h\|_{\tilde\patch}\|\nabla (z-z_H)\|_{\tilde\patch}.
  \end{aligned}
  \]
  Hence, taking $z_H=I_H z\in V_{H,0}(\tilde\patch)$ in combination with a corresponding error estimate and the estimate of the first part yields
  \[
  \big|(v_h,\psi)_{\tilde\patch}\big|
  \le cH\|\nabla (u_h-u_H) \|_{\tilde\patch} \|\nabla^2 z\|_{\tilde\patch}
  \le cH\|\nabla (u_h-u_H) \|_{\tilde\patch} \|\psi\|_{\tilde\patch},
  \]
  where we used \eqref{eq:H2_aprior_patch} in the last step. Taking the supremum over $\psi\in L^2(\tilde\patch)$ gives the second estimate of the assertion.
\end{proof}

After these preparations we conclude with the proof of the main theorem.
\begin{proof}[Proof of Theorem \ref{thrm2}]
  As in the single-patch case we first introduce the fine mesh solution $u_h\in V_h$ to the right hand side $f$ and a corresponding error estimate yielding
  \begin{equation}\label{e2:1}
    \|\nabla (u-u_\NN) \| \le
    c h^r \|f\|_{r-1}+
    \|\nabla (u_h-u_\NN) \|.  
  \end{equation}
  As $u_\NN=u_H+\sum P_\patch w_\NN^{\patch}$ is composed of local updates, we, from here on, consider the local contributions $\|\nabla (u_h-u_\NN)\|_\patch$.
  On each patch $\patch$ we introduce the following local solutions to local problems given by Definition~\ref{def:local}
  \begin{equation}\label{locproblems}
    \begin{aligned}
      v_h&\in V_{h,0}(\tilde\patch)&\quad (\nabla (u_H+v_h),\nabla\phi_h)_{\tilde\patch} &= (f,\phi_h)_{\tilde\patch}&\quad\forall \phi_h&\in V_{h,0}(\tilde\patch),\\
      v_h^{i}&\in V_{h,0}(\tilde\patch)& (\nabla (u_H^i+v_h^{i}),\nabla\phi_h)_{\tilde\patch} &= (f_i,\phi_h)_{\tilde\patch}&\forall \phi_h&\in V_{h,0}(\tilde\patch). 
    \end{aligned}
  \end{equation}
  The index $i$ will refer to an element of the training data. The first local problem in \eqref{locproblems} is to estimate the local error between $u_h$ and $u_H^h:=u_H+v_h$, whereas the second local problem estimates the local error between the training data $u_h^i$ and $u_H^{h,i}:=u_H^i + v_h^{i}$.  We split the error as
  \begin{multline}\label{loc2}
    \|\nabla (u_h-u_\NN)\|_\patch \le
    \|\nabla (u_h-u_H^h)\|_\patch
    +\|\nabla (u_H^h-u_H^{h,i})\|_\patch\\
    +\|\nabla (u_H^{h,i}-u_h^i)\|_\patch
    +\|\nabla (u_h^i-u_\NN)\|_\patch.
  \end{multline}
  In principal, the first and third term can be estimated by~\cite[Theorem 5]{XuZhou1999}, where local finite element algorithms are analyzed that decompose the solution $u_h$ into a global coarse solution $u_H\in V_H\subset V_h$ and into a local fine mesh solution. However, as $u_h-u_H^h$ is discretely harmonic on $\tilde\patch$, we can just use~\cite[Lemma 3.3]{DemlowGuzmanSchatz2011}, which is the main ingredient for the proof of local finite element error estimates, introduce $\pm u_H$ and use Lemma~\ref{lemma:local} to get
  \begin{equation}\label{loc3} 
    \begin{aligned}
      \|\nabla (u_h-u_H^h)\|_\patch & \le c d(\patch,\tilde\patch)^{-1} \|u_h-u_H^h\|_{\tilde\patch}\\
      &\le c h_{\patch}^{-1}\big(\|u_h-u_H\|_{\tilde\patch}
      + H\|\nabla (u_h-u_H)\|_{\tilde\patch}\big).
    \end{aligned}
  \end{equation}
  Note that in the last step we used $d(\patch,\tilde\patch)\sim h_{\patch}$, which is valid by construction.
  Likewise, for $u_h^i-u_H^{h,i}$ in~\eqref{loc2} we get
  \begin{equation}\label{loc5a}
      \|\nabla (u_H^{h,i}-u_h^i)\|_\patch \le c h_{\patch}^{-1}\big(\|u_h^i-u_H^i\|_{\tilde\patch}
      + H\|\nabla (u_h^i-u_H^i)\|_{\tilde\patch}\big).
  \end{equation}
  Next, considering the definition of~\eqref{locproblems}, the second term in~\eqref{loc2} can be estimated as an error that measures the richness of the training data, i.e., we obtain
  \[
    \begin{aligned}
        \|\nabla (u_H^h-u_H^{h,i})\|_{\patch}^2&\le \|\nabla (u_H^h-u_H^{h,i})\|_{\tilde \patch}^2\\
        %&=(\nabla(u_H^h-u_H^{h,i}),\nabla(u_H-u_H^i))_{\tilde\patch}+(\nabla(u_H^h-u_H^{h,i}),\nabla(v_h-v_h^i))_{\tilde\patch}\\
        &=(\nabla(u_H^h-u_H^{h,i}),\nabla(u_H-u_H^i))_{\tilde\patch}+(f-f_i,v_h-v_h^i)_{\tilde\patch}\\
        %&\le \|\nabla(u_H^h-u_H^{h,i})\|_{\tilde\patch}\|\nabla(u_H-u_H^i)\|_{\tilde\patch}+c\|f-f_i\|_{-1,\tilde\patch}\|\nabla(v_h-v_h^i)\|_{\tilde\patch}\\
        &\le \|\nabla(u_H^h-u_H^{h,i})\|_{\tilde\patch}\|\nabla(u_H-u_H^i)\|_{\tilde\patch}\\
        &\quad+c\|f-f_i\|_{-1,\tilde\patch}\Big(\|\nabla(u_H^h-u_H^{h,i})\|_{\tilde\patch}+\|\nabla(u_H-u_H^i)\|_{\tilde\patch}\Big),
    \end{aligned}
  \]
  where we inserted $\pm (u_H-u_H^i)$ as an intermediate function.
  As a consequence of Young's inequality, we get
  \begin{equation}\label{loc5}
    \|\nabla (u_H^h-u_H^{h,i})\|_{\patch}\le c\left(\|\nabla (u_H-u_H^i)\|_{\tilde \patch} + \|f-f_i\|_{-1,\tilde\patch}\right).
  \end{equation}
  Note that the constant $c$ does not depend in $h_{\patch}$.
  The last term in~\eqref{loc2} is split into
  \begin{equation}\label{loc6}
    \|\nabla (u_h^i-u_\NN)\|_\patch\le
    \|\nabla (u_h^i-u_\NN^i)\|_\patch+\|\nabla (u_\NN^i-u_\NN)\|_\patch,
  \end{equation}
  where $u_\NN^i:= u_H^i+ \sum_{\patch}P_{\patch}w_\NN^{\patch, i}$ is the local network approximation for the training data elements $u_H^i$ and $f_i$. By \eqref{eq:norms2} and the assumption that the loss is reduced to the order of $\epsilon^2$ for each patch for each training problem, we get a bound for the first term, i.e.,
  \begin{equation}\label{loc61}
    \|\nabla (u_h^i-u_\NN^i)\|_\patch 
    \le c h^{\frac{d}{2}-1} \|u_h^{i}-u_\NN^{i}\|_{\ell^2(\patch )} \le c h^{\frac{d}{2}-1} \epsilon,
  \end{equation}
  which depends on the expressivity of the neural network and the optimization error.
  The second term of~\eqref{loc6} is estimated as
  \begin{equation}\label{loc7}
    \|\nabla (u_\NN^i-u_\NN)\|_\patch
    \le \|\nabla (u_H^i-u_H)\|_\patch+ c\|\nabla (w_\NN^{\patch, i}-w_\NN^{ \patch})\|_\patch
  \end{equation}
  and again consists of the data error $\|\nabla(u_H-u_H^i)\|_{\patch}$ and, finally, the local network generalization error.
  We combine~\eqref{loc2}-\eqref{loc7} and sum over all patches to get
  \[
    \begin{aligned}
    \|\nabla (u_h-u_\NN)\|^2 &\le c\Big[
      h_{\patch}^{-2}\left(\|u_h-u_H\|^2+H^2\|\nabla(u_h-u_H)\|^2\right)\\
      &\quad+\sum_{\patch\in U_\patch}
      \min\limits_{i \in \overline{1,N_T}}\Big\{h_{\patch}^{-2}\left(\|u_h^i-u_H^i\|_{\tilde\patch}^2+H^2\|\nabla(u_h^i-u_H^i)\|_{\tilde\patch}^2\right)\\
      &\quad+\|\nabla (u_H-u_H^i)\|_{\tilde \patch}^2 + \|f-f_i\|_{-1,\tilde\patch}^2
      +\|\nabla (w_\NN^{\patch, i}-w_\NN^{\patch})\|_\patch^2\Big\}+N_Ph^{d-2} \epsilon^2
      \Big]
    \end{aligned}
  \]
  The first two terms can be estimated by means of \eqref{femerror} after having introduced $\pm u$ as intermediate functions. The network update term can be bounded in the same way as in \autoref{thrm1} by using Lemma \ref{lemma:network_stability}. Finally, using $h\le H$ and taking the square root then yields
\begin{align*}
    \norm{\nabla  (u_h - u_\NN)}
    &\leq c \Bigg[
        H^{r+1} h_\patch^{-1}  \norm{f}_{r-1}
        +\Big(\sum\limits_{\patch \in U_\patch} \min\limits_{i \in \overline{1,N_T}} \Big\{
            h_\patch^{-2} \Big( \norm{u_h^i - u_H^i}_{\tilde{\patch}}^2\\
            &\quad+ H^2 \norm{\nabla (u_h^i - u_H^i)}_{\tilde{\patch}}^2\Big)
            + \norm{\nabla(u_H - u_H^i)}_{\tilde{\patch}}^2 + \norm{f-f_i}_{-1, \patch}^2 \\
            & \quad + L(\NN)^2h^{-2} \Big(
                \norm{u_H - u_H^i}_\patch^2 + \norm{I_h(f-f_i)}_\patch^2
            \Big)
        \Big\} \Big)^{\frac{1}{2}} + N_\patch^\frac12 h^{\frac{d}{2}-1}\epsilon
    \Bigg],
\end{align*}
hence giving the final result.  
\end{proof}

\begin{remark}
The error estimates from Theorems~\ref{thrm1} and~\ref{thrm2} show a balanced error. In particular, they prove the practical usefulness of the method: the total error is given as a balance of richness of the training data, i.e., number and distribution of the training elements, and quality of the training data, i.e., the resolution of the training data. Both can be controlled. In particular, here $u_H^i, u_h^i, f_i$ are training data and $u_H$ and $f$ are known beforehand and hence all terms involving only these functions can also be computed a priori. Due to that we do not modify them any further and leave them as they are in the estimate. In addition, there is the network and optimization error, which in principle can be influenced by the architecture, depth and width of the network. The practical efficiency of the method will always depend on the specific example and especially on how much the effort of offline phase and online phases differ. For typical problems in fluid mechanics, the 3d flow around obstacles, we observed substantial increases in efficiency for relevant generalizations of the training data~\cite{Margenberg2023}. 

Numerically we could not fully identify proper relation between the mesh sizes $h$ and $H$ as well as the patch size $h_P$. The estimates suggest that it must hold $H/h_{\cal P}=o(1)$ such that the method gives real benefit, e.g. by choosing $h_{\cal P}=\sqrt{H}$. Our numerical experiments however show optimal performance for $H/h_{\cal P}=O(1)$. In future work we will aim at improving the theoretical results, e.g. by exploiting super-convergence estimates in the locally structured meshes.
\end{remark}

\begin{remark}[Extension to training domains]\label{remark:uniformmeshes}
Instead of using the domain $\Omega$ for generating the training data as well as for calculating the hybrid solution, one could also use a training domain $\Omega^{tr}$ instead, which differs from the application domain $\Omega$. This would give us more flexibility during the whole training process.

At first, we could consider the following situation:
We denote by $\Omega^{tr}$ the domain on which the training data is generated, 
and by $\Omega^{tr}_h$, $\Omega^{tr}_H$, and $\Omega^{tr}_{h_\patch}$ the corresponding meshes.
We assume that the set of patches
\[
  U_\patch = \{\patch_1,\dots,\patch_{N_\patch}\}
\]
coincides with the cells of $\Omega_{h_\patch}$, i.e., $U_\patch = \Omega_{h_\patch}$.
Similarly, the set of training patches
\[
  U_\patch^{tr} = \{\patch^{tr}_1,\dots,\patch^{tr}_{N_\patch^{tr}}\}
\]
coincides with the cells of $\Omega_{h_\patch}^{tr}$, i.e., 
$U_\patch^{tr} = \Omega_{h_\patch}^{tr}$.
Moreover, we assume that $U_\patch$ and $U_\patch^{tr}$ are compatible in the following sense: 
for each $\patch \in U_\patch$ there exists a training patch 
$\patch^{tr} \in U_\patch^{tr}$ that is a translation and/or rotation of $\patch$. An illustration is given in Fig.~\ref{fig:meshes}.
\begin{figure}[t]
  \begin{center}
    \includegraphics[width=0.8\textwidth]{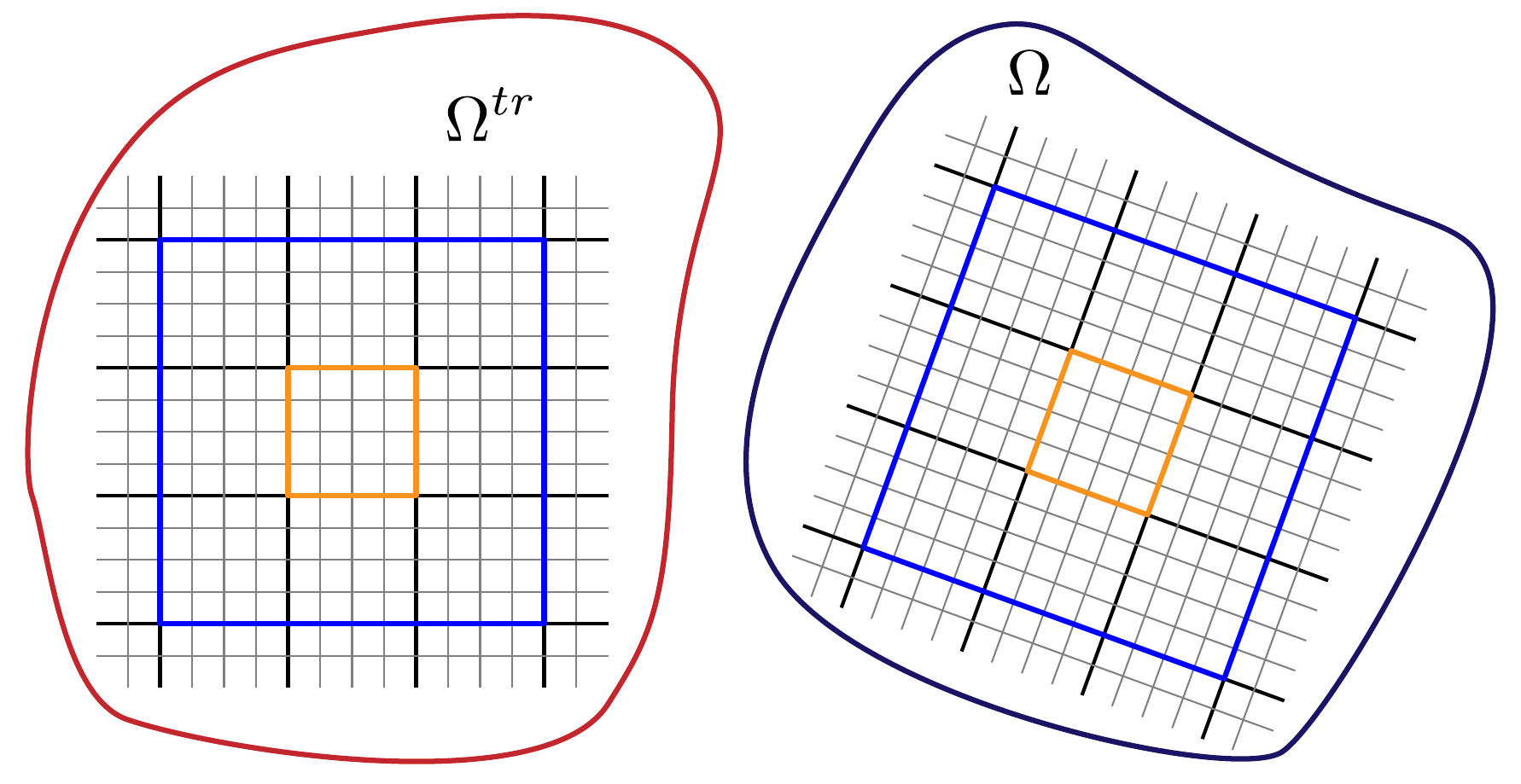}
  \end{center}
  \caption{Training and application domains $\Omega^{tr}$ and $\Omega$ as well as the corresponding meshes $\Omega^{tr}_h$ and $\Omega_h$ can differ. Both however must be split into the same kind of patches. A patch $\patch$ and the extended patch $\tilde\patch$ is marked in both domains in orange and blue, respectively.}
  \label{fig:meshes}
\end{figure}

This assumption is restrictive and only allows fully uniform meshes.
Hence, it gives little freedom in the generalization of the domain $\Omega^{tr}$. Basically, all domains have to be put together by blocks of patches that are found in the training data. The path for a further generalization would be by means of introducing a parametric setup: data is still kept local on the patches $\patch$, but the training of the network and the application of the network is always by means of transformation to a reference patch $\patch_r$. In setting up the neural network approach, the training data then must include a sufficient variety of different patch geometries and sizes that are also found in the application domain $\Omega$. This approach would also be a first step towards an application to locally refined finite element methods. 
\end{remark}
%

% %%%%%%%%%%%%%%%%%%%%%%%%%%%%%%%%%%%%%%%%%%%%%%%%%%

\section{Numerical experiments}\label{sec:num}

In the following paragraphs we will document different numerical simulations in order to explore the performance of the hybrid finite element neural network approach. In particular, we will study the sensitivity of the approximation properties on the various aspects like network size, variety of training data or data preprocessing techniques. Then, we consider a second test case with a different geometry to explore the generalization capacity of the approach.

\subsection{Configuration of the experimental setup}\label{subsec:setup}
\begin{figure}[htbp]
    \centering
    \begin{tabular}{ccc}
        \includegraphics[width=0.3\linewidth]{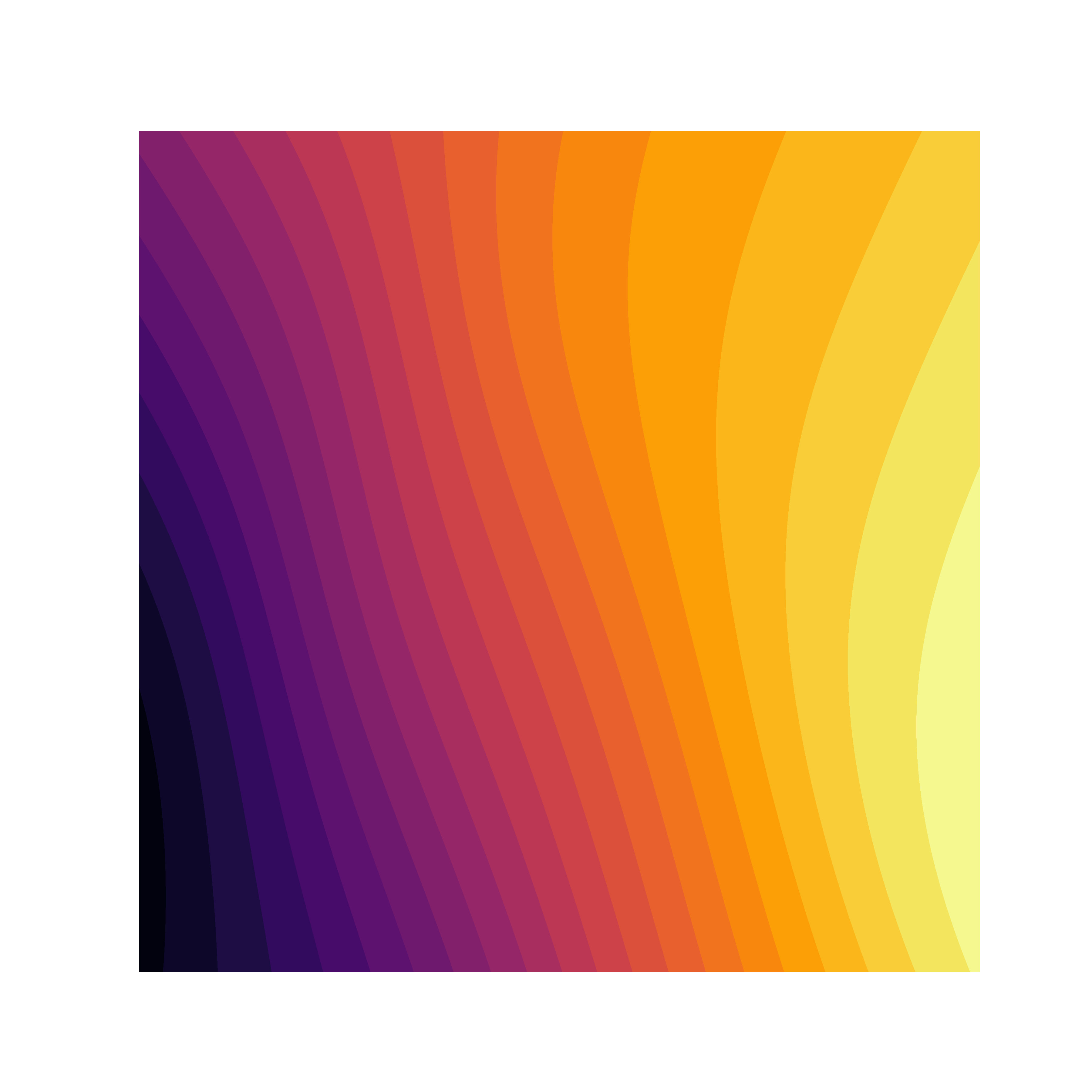} & 
        \includegraphics[width=0.3\linewidth]{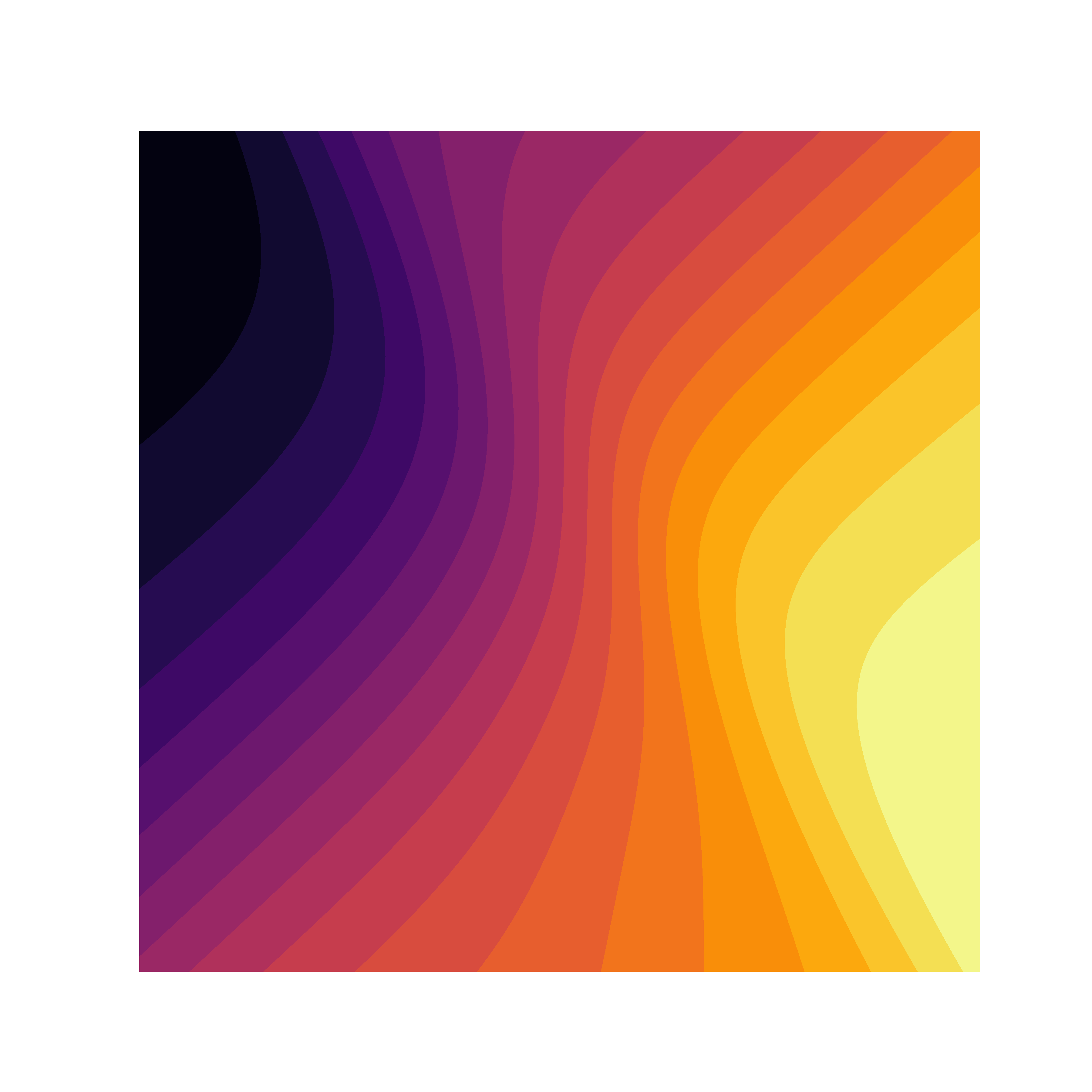} & 
        \includegraphics[width=0.3\linewidth]{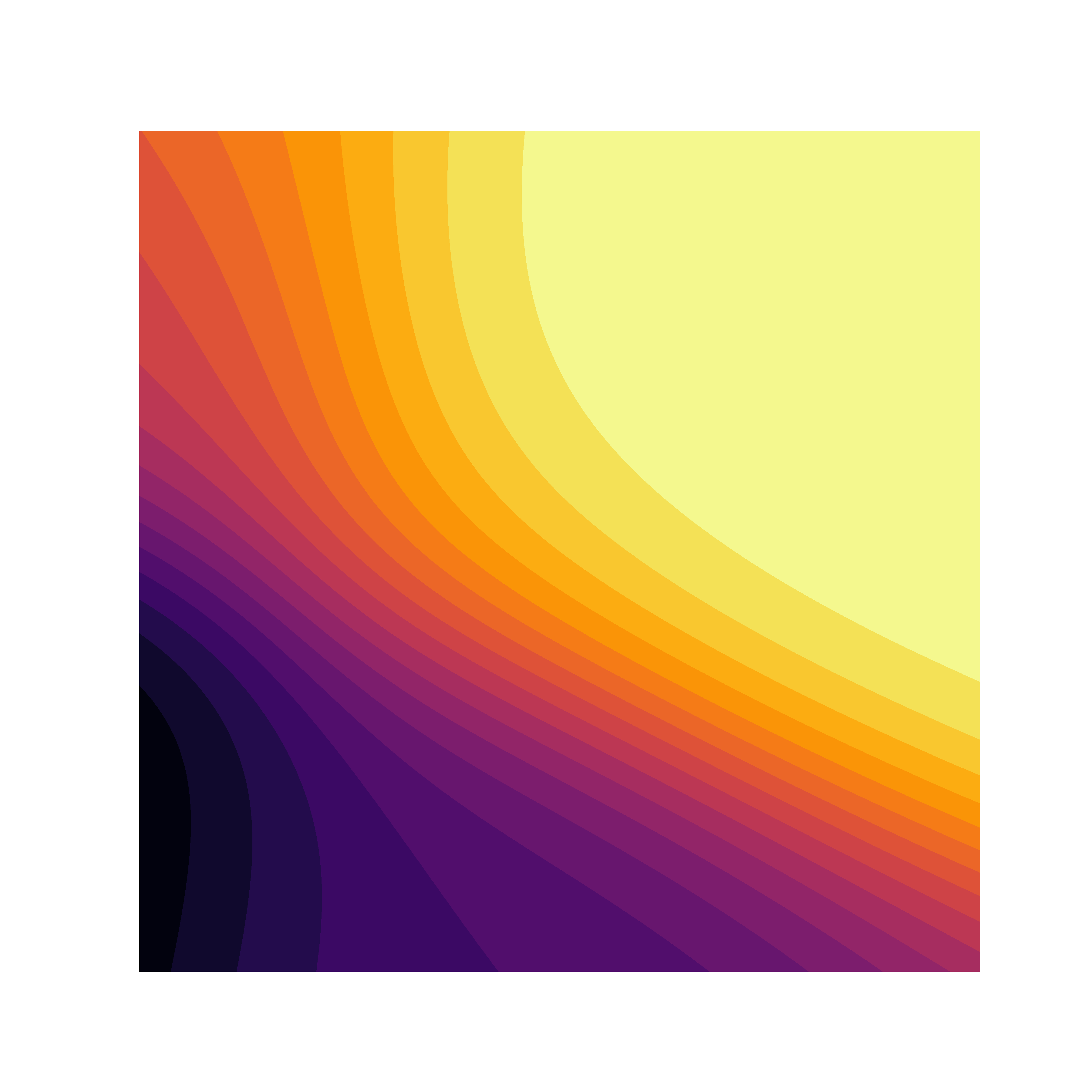} \\
        \includegraphics[width=0.3\linewidth]{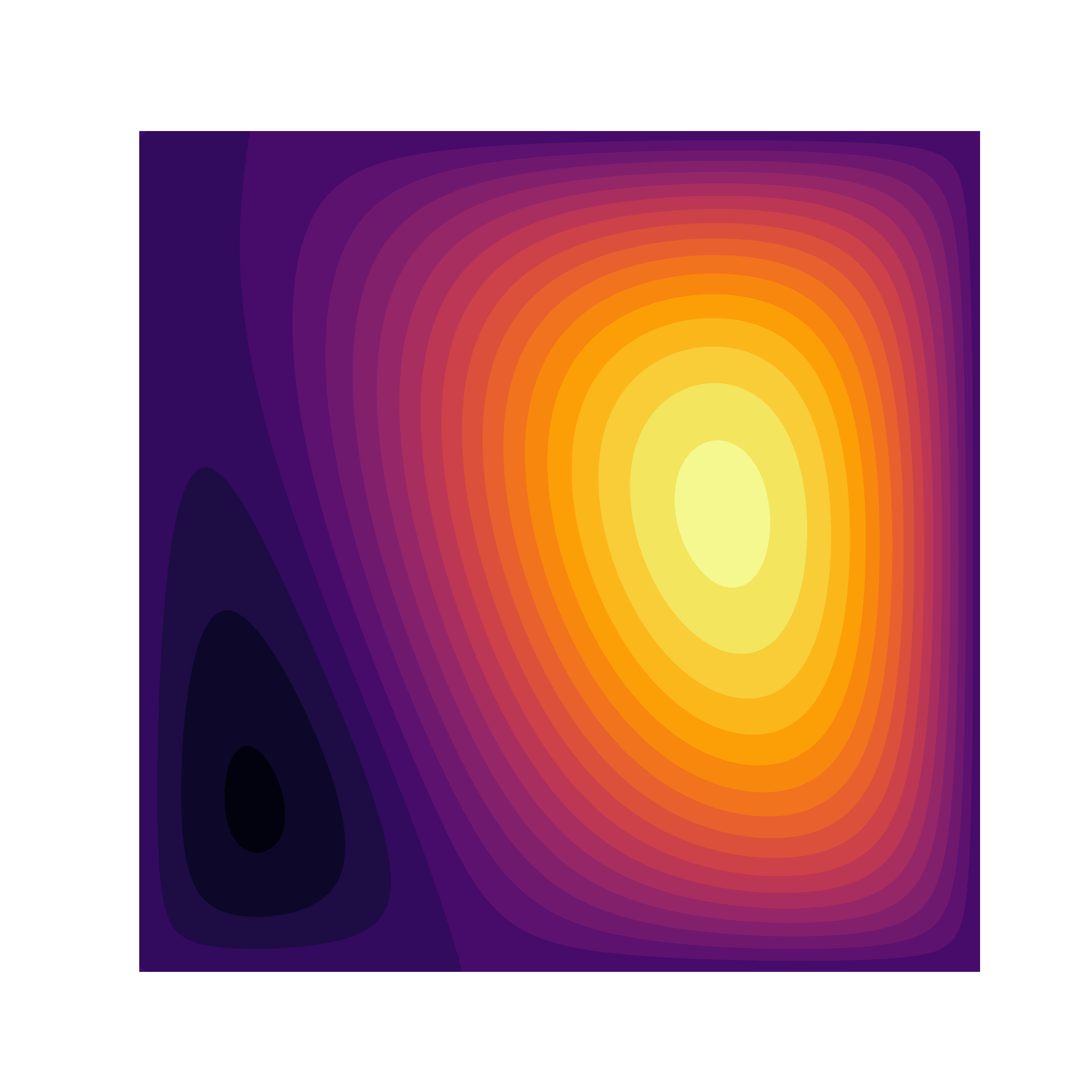} & 
        \includegraphics[width=0.3\linewidth]{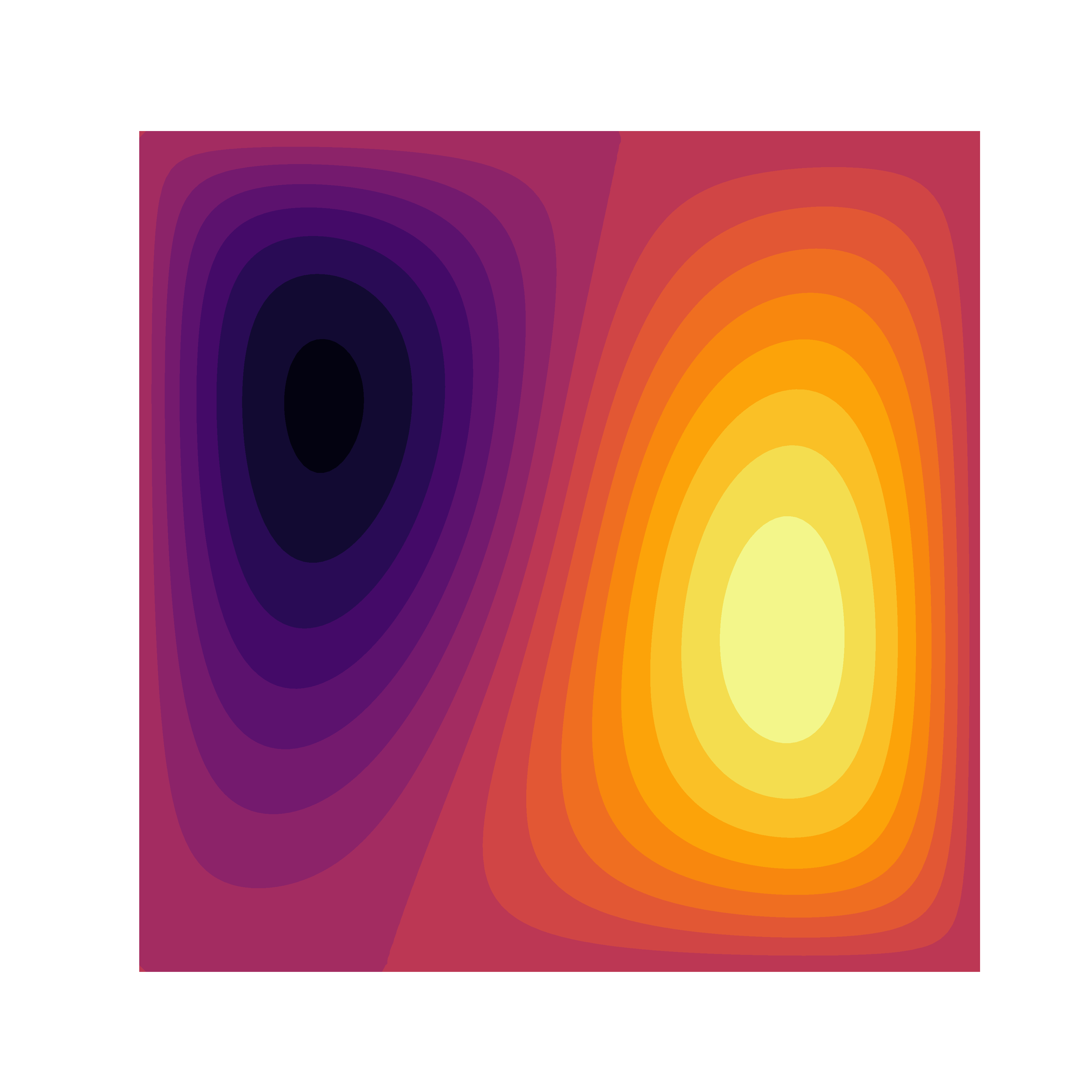} & 
        \includegraphics[width=0.3\linewidth]{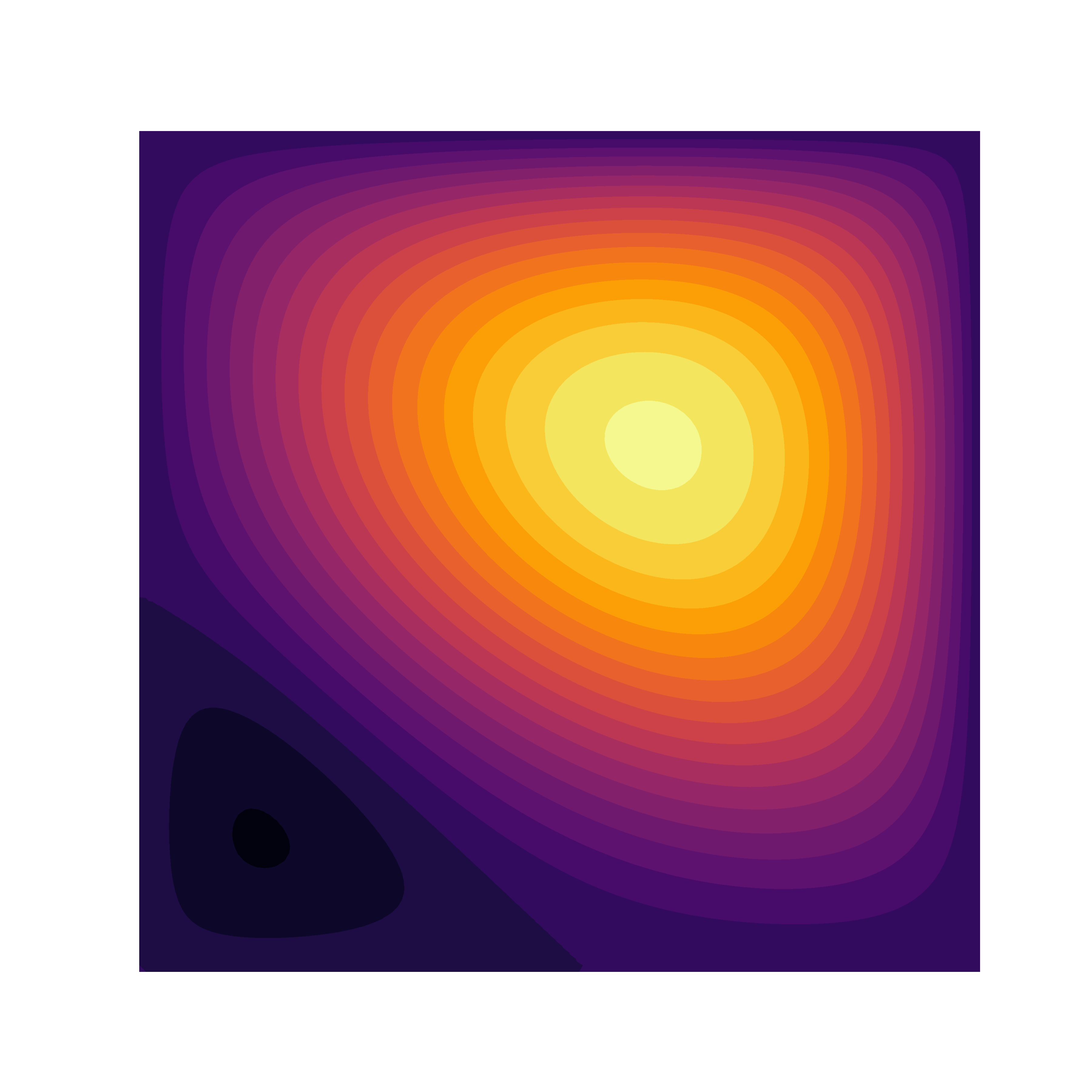} \\
    \end{tabular}
\caption{
Top: example source terms from ${\cal F}$. Bottom: corresponding solutions}
\label{fig:example_rhs_and_solutions}
\end{figure}

We start by describing our experimental setup. Let $\Omega=(0,1)^2$ be the unit square. We consider the two-dimensional Poisson equation with homogeneous Dirichlet boundary conditions
\begin{equation}
  -\Delta u = f\text{ in }\Omega,\quad
  u = 0\text{ on }\partial\Omega,
\end{equation}
where $f$ is an element of the following set
\begin{equation}\label{Fset}
\begin{aligned}
  {\cal F} := \Big\{\, 
  f(x) &= \sum_{i=1}^{n} \alpha_i \sigma\!\big(s_i ( W_i x + b_i )\big)
  : W_i \in \mathbb{R}^{2 \times 2},\; b_i \in \mathbb{R}^2, \\[3pt]
  &\alpha_i \in [0,1],\; \sum\nolimits_{i=1}^n \alpha_i = 1,\;
  s_i \in [1,20]
  \,\Big\}
\end{aligned}
\end{equation}
where $W_i$  and $b_i$ define a hyperplane going through some points in $\Omega$. The source terms for both training and testing are generated by picking at random samples from this set. In \autoref{fig:example_rhs_and_solutions} one can see examples of such source terms and corresponding solutions.

% \begin{figure}
% \begin{center}
% \end{center}
% \caption{Solution examples}
% \end{figure}

\subsection{Neural network setup and optimization}

In this section we consider a  multilayer perceptron as described previously in Section~\ref{sec:hybrid}, see Definition~\ref{def:MLP}. We consider networks with $L=1$ hidden layer of size 64, if not stated otherwise. The total number of trainable parameters depends on the size of the input and output vectors, see Table \ref{table:params}. In order to predict the correction for the fine level, the network receives $4$ nodal values of the coarse solution and $(2^k+1)^2$ nodal values of the source term. The output vector contains $(2^k+1)^2$ nodal values which predicts the difference between fine and coarse finite element solutions. In general, for multilayer perceptrons, the number of inputs and outputs, as well as the total number of trainabale parameters, can be computed as follows

\begin{equation}
\begin{aligned}
    \#inputs    := & \left(2^{\left(CL-PL\right)}+1\right)^2 +  \left(2^{\left(FL-PL\right)}+1\right)^2 \\
    \#outputs   := & \left(2^{\left(FL-PL\right)}+1\right)^2 \\
    \#trainable_parameters:= & (\#inputs + 1) \cdot \#neurons \\
      & + (\#layers-1) \cdot \#neurons \cdot (\#neurons + 1) \\
      & + (\#neurons+1) \cdot \#outputs \\
\end{aligned}
\end{equation} 

where $PL$, $CL$ and $FL$ denote patch, coarse and fine levels, respectively.

\begin{table}[t]
  \centering
  \begin{tabular}{c|c|c|c|c|c}
    \toprule
    H/h   & $N_0$    & $N_1$ & $N_2$ & \# parameters &  \#DOF\\ 
    \midrule
    $2^1$ &  4 +  25 &  64                & 25 &     1\,536 & 1089\\
    $2^2$ &  4 +  81 &  64                & 81 &     3\,584 & 4225\\
    $2^3$ &  4 + 289 &  64                & 289  &  10\,752 & 16641\\
    \bottomrule
  \end{tabular}
  \caption{The number of neurons $N_l$ at each layer $l \in \{0,\dots,2\}$, the total number of trainable parameters (weights and biases) depending on the jump level and the corresponding total number of degrees of freedom (cl=3)}
  \label{table:params}  
\end{table}

\subsubsection{Generation of training data}

For creating training and test data sets we randomly generate finite sets of  functions taken from  ${\cal F}$ as given in~\eqref{Fset}.
% To generate such a function it suffices to uniformly sample $d \cdot n$ points from the domain for computation of $W_i$ and $b_i$, $n$ numbers from $[0,1]$ for $\alpha_i$ and $n$ real numbers from $[1,10]$ for $s_i$.

The computations are performed for a coarse step size $H$ and for the fine mesh size $h=H/8$. Hereby, we generate sets of training and test triples $(f_i,u_H^i,u_h^i)\in {\cal F}\times V_H\times V_h$ for $i=1,\dots,N_T$. In our experiments, we set $H=2^{-4}$ which corresponds to $256$ cells.

Due to the local neural network setup, each triple $(f_i,u_H^i,u_h^i)$ generates a larger number of training patches $N_\patch$, one for each patch of the triangulation. The number of patches are always the same as the number of elements in the coarse mesh. Hence, $H=h_\patch$ and we have $N_\patch=256$.

A second test data set of the same size is additionally chosen from ${\cal F}$. 

\FloatBarrier

\subsubsection{Optimization}

For the training of the network we use the Adam optimizer~\cite{kingma2017adam}. The loss function is the mean square error as defined in \eqref{minimize} averaged by the batch size. We train the network for 10000 epochs with constant learning rate $10^{-3}$.

% %%%%%%%%%%%%%%%%%%%%%%%%
\subsection{Accuracy of the hybrid finite element neural network solver}

In this first test case, we analyze the accuracy of the hybrid solver. For this purpose, we consider the mean error over all test data. Starting from the coarse solution with grid size $H=2^{-4}$ and $h_\patch = H$, the hybrid method is enriched with neural networks. In doing so, we predict 3 grid levels, thus trying to achieve accuracy of $h=H/8$.  For evaluating quality of solutions we use $\ell_2$-norm as defined in \eqref{eq:l2_norm}, but additionally averaged over degrees of freedom, so that it is consistent over different meshes, which makes it the same as averaged root mean square error
\begin{equation*}
    RMSE(x, y) := \frac{1}{\sqrt{N}} \norm{x-y}_2,
\end{equation*}
where $N$ is a number of degrees of freedom, and $x, y\in\mathbb{R}^N$ are vectors of nodal values.

In this first test, we only want to investigate the approximation ability of the neural network and we use an excessive number of training data to keep the data error term 
\[
\min_i \big\{\|f-f_i\| + \|\nabla (u_H - u_H^i)\|\big\}
\]
small. Thus, from~\eqref{aprioriestimate}, considering linear finite elements $r=1$, it remains
% \highlight{
\[
\|\nabla (u - u_\NN)\| = {\cal O}\big(h+H^2 h_{\patch}^{-1}\big).
\]
% }

Since we do not know the analytical solution we use a reference solution $u_{\text{ref}}$ instead, which is calculated on a mesh with element size  $h_{\text{ref}}:=h/2$. The numerical results shown in Table~\ref{tab:n1} are superior to the theoretical estimate. When estimating one or two mesh levels, the full fine mesh accuracy is reached and the discrepancy to the fine mesh error is still very small for estimating three mesh levels. The negative impact of the patch size $h_{\cal P}^{-1}$ does not strongly show up. 
% {\color{red}VLAD: please write down what exactly MSE is.} {\color{orange}All errors are indicated as the mean squared error
% \[
% MSE = ...
% \]}
%
\begin{table}
  \begin{center}
    \begin{minipage}{0.45\textwidth}
      \begin{tabular}{cccc}
        % data/error_scaling_wrt_fine_level.csv 
        \toprule
        $h$ & $\|u-u_h\|$& $\|u-u_\NN\|$\\
        \midrule
        $H/2$ & 1.51e-05 & 1.51e-05 \\
        $H/4$ & 3.82e-06 & 3.77e-06 \\
        $H/8$ & 9.60e-07 & 1.11e-06 \\
        \bottomrule 
      \end{tabular}
    \end{minipage}
    \begin{minipage}{0.44\textwidth}
      ~\hspace{-4cm}\vspace{-0.5cm}
      \begin{tikzpicture}[scale=1]
        \begin{axis}[
            width=0.8\textwidth, 
            height=0.6\textwidth,
            legend style ={at={(-0.3,-0.1)},anchor=north east},
            legend columns=-1,
            xlabel={Relative mesh size $h/H$},
            ylabel={Error},
            xmode=log,
            log basis x={2},
            ymode=log,
            log basis y={2},
            mark options={solid, scale=1.6}
          ]
          \addplot[very thick, dashed, blue,  mark=Mercedes star] coordinates {(0.5,1.51e-05) (0.25, 3.82e-06) (0.125,9.60e-07)};
          \addplot[very thick, olive, mark=Mercedes star flipped] coordinates {(0.5,1.51e-5) (0.25, 3.77e-06) (0.125,1.11e-06)};
          \legend{
            Reference  $\|u-u_h\|$,
            Hybrid  $\|u-u_\NN\|$
          }
        \end{axis}
      \end{tikzpicture}
    \end{minipage}
  \end{center}
  \vspace{0.5cm}
  \caption{Performance of the hybrid solution for three different refinement levels. We show the average error on the test data. The training data set ${\cal F}$ is sufficiently big such that the data error is negligible.}
  \label{tab:n1}
\end{table}

% %%%%%%%%%%%%%%%%%%%%%%%%

\subsubsection{Dependency on the training data set}
\begin{figure}[t]
\begin{tikzpicture}

\begin{groupplot}[
    group style={
        group size=3 by 1,
        horizontal sep=.7cm,
        vertical sep=1cm,
    },
    width=4.5cm,
    height=3cm,
    xmode=log,
    log basis x={2},
    ymode=log,
    log basis y={2},
]

\nextgroupplot[
    legend to name=commonlegend,
    legend columns=-1,
    legend style={draw=black, fill=white},
    legend image post style={draw=none},
    xlabel=$N_T$,
    ylabel=Error
]

\addplot[
    mark=square*, 
    mark options={solid, scale=1.4}, 
    red, 
    dotted, 
    very thick
] table [
    x=n_rhss, 
    y=uff-uc (test), 
    col sep=comma
] {data/data_scaling.csv};
\addlegendentry{$\|u_H-u\|$}
\addplot[
    mark=Mercedes star, 
    cyan,  
    very thick
] table [
    x=n_rhss, 
    y=uff-uf (test), 
    col sep=comma
] {data/data_scaling.csv};
\addlegendentry{$\|u_{H/8} - u\|$}
\addplot[
    mark=Mercedes star flipped, 
    brown, 
    dashed, 
    very thick
] table [
    x=n_rhss, 
    y=uff-un (test), 
    col sep=comma
] {data/data_scaling.csv};
\addlegendentry{$\| u_\NN - u \|$}

\nextgroupplot[
    yticklabels=\empty,
    xlabel=Number of layers
]
\addplot[
    mark=square*, 
    mark options={solid, scale=1.4}, 
    red, 
    dotted, 
    very thick
] table [
    x=num_layers, 
    y=uff-uc (test), 
    col sep=comma
] {data/num_layers_scaling.csv};
\addplot[
    mark=Mercedes star, 
    cyan, 
    very thick
] table [
    x=num_layers, 
    y=uff-uf (test), 
    col sep=comma
] {data/num_layers_scaling.csv};
\addplot[
    mark=Mercedes star flipped, 
    brown, 
    dashed, 
    very thick
] table [
    x=num_layers, 
    y=uff-un (test), 
    col sep=comma
] {data/num_layers_scaling.csv};

\nextgroupplot[
    yticklabels=\empty,
    xlabel=Number of neurons
]
\addplot[
    mark=square*, 
    mark options={solid, scale=1.4}, 
    red, 
    dotted, 
    very thick
] table [
    x=num_neurons, 
    y=uff-uc (test), 
    col sep=comma
] {data/num_neurons_scaling.csv}; 
\addplot[
    mark=Mercedes star, 
    cyan,  
    very thick
] table [
    x=num_neurons, 
    y=uff-uf (test), 
    col sep=comma
] {data/num_neurons_scaling.csv};
\addplot[
    mark=Mercedes star flipped, 
    brown, 
    dashed, 
    very thick
] table [
    x=num_neurons, 
    y=uff-un (test), 
    col sep=comma
] {data/num_neurons_scaling.csv};

\end{groupplot}

% Place the shared legend below the group
\node at ($(group c2r1.south) + (0,-1.5cm)$) {\ref{commonlegend}};

\end{tikzpicture}
\caption{Dependency of the prediction quality $\|u-u_\NN\|$ (for the test data) on the size of the training data set, number of layers and number of neurons}
\label{fig:performance_data_amount}
\end{figure}

Next, in the left plot of Figure~\ref{fig:performance_data_amount} we show the approximation error of the hybrid simulation for jumping three mesh levels, i.e. $h=H/8$, depending on the size of the training data set $N_T$. The results indicate that the corresponding fine mesh accuracy is reached, however, only for an increased amount of training data. This is in accordance to the error estimate~\eqref{aprioriestimate}, which shows that the error term $h^r + H^{r+1} h_{\patch}^{-1}$ must be balanced with the data error $\min\limits_i \|f-f_i\|$ and $\min\limits_i \|\nabla (u_H-u_H^i)\|$. Furthermore, Figure~\ref{fig:performance_data_amount} suggests the scaling $\min\limits_i \|f-f_i\| = {\cal O}(N_T^{-\frac{1}{2}})$.

% %%%%%%%%%%%%%%%%%%%%%%%%
\subsubsection{Error depending on the network complexity}

% A goal for opimization is to add bounds on $c_W$ to the loss function with the goal to keep this essential constant small. In the next subsection we will present the obtained results. 

Next we study the dependency of the error on the complexity of the network. Again we consider the case  $h=H/8$. We used $N_T=2^{9}$ training problems for this experiment. First we choose the number of hidden layers of the perceptron between 1, 2 or 4 with 64 neurons in each layer. 

The results can be seen in the middle plot of \autoref{fig:performance_data_amount}. We observe that most of the time using one hidden layer seems to be sufficient.

Finally, we fix the number of layers to be one and vary the number of neurons per layer from 64 to 512. The results can be seen in the right plot of  \autoref{fig:performance_data_amount}. Here we observe that the error behaves similarly to the experiment where we varied the amount of training data - the error decreases as the number of neurons increases. We confirm that 64 is most likely indeed the optimal choice in the current setting. 

\subsection{Role of data preprocessing}

Finally, we justify our specific choice of data preprocessing, so-called standardization, which we employed in every experiment presented in this paper. The data preprocessing has a large effect on the optimization of the network and its approximation quality. In order to illustrate this, we consider fixed $h = H/8$ and we vary the amount of training data from $2^{6}$ to $2^{9}$. 
As for preprocessing methods, we considered no preprocessing and and standardization. These methods were applied to each input and output feature separately. Standardization transforms the data to have mean value of $0$ and variance of $1$, i.e.
\begin{equation}
  \tilde{X}_{std} = \frac{X - \operatorname{Mean}(X)}{\operatorname{Var}(X)}.
\end{equation}

\begin{figure}
\begin{tikzpicture}
\begin{groupplot}[
    group style={
        group size=2 by 1,
        horizontal sep=.5cm,
        vertical sep=1cm,
    },
    width=6cm,
    height=3cm,
    ymode=log,
    log basis y={2},
    ymax=0.00009
]

\nextgroupplot[
    legend to name=commonlegend2,
    legend columns=-1,
    legend style={draw=black, fill=white},
    legend image post style={draw=none},
    xmode=log,
    log basis x={2},
]

\addplot[very thick, red, mark=square*, mark options={solid,scale=1.2}] table [x=n_rhss, y=uff-uc (test), col sep=comma, solid] {data/preprocessing.csv};
\addlegendentry{$\|u_H-u\|$}
\addplot[very thick, blue, mark=triangle*] table [x=n_rhss, y=uff-uf (test), col sep=comma, solid] {data/preprocessing.csv};
\addlegendentry{$\|u_h-u\|$}
\addplot[very thick, brown, mark=Mercedes star] table [x=n_rhss, y=No preprocessing (test), col sep=comma] {data/preprocessing.csv};
\addlegendentry{$\|u_\NN-u\|$ (no preprocessing)}
\addplot[very thick, orange, mark=Mercedes star] table [x=n_rhss, y=Standard scaler (test), col sep=comma] {data/preprocessing.csv}; 
\addlegendentry{$\|u_\NN-u\|$ (standardization)}
		  
\nextgroupplot[yticklabels=\empty]
		  
\addplot[very thick, red] table [x=epoch, y=uff-uc (test), col sep=comma, solid] {data/error_during_training_wrt_preprocessing.csv}; 
\addplot[very thick, blue] table [x=epoch, y=uff-uf (test), col sep=comma, solid] {data/error_during_training_wrt_preprocessing.csv};
\addplot[very thick, brown] table [x=epoch, y=uff-un (test no preprocessing), col sep=comma] {data/error_during_training_wrt_preprocessing.csv};
\addplot[very thick, orange] table [x=epoch, y=uff-un (test standard scaler), col sep=comma]   {data/error_during_training_wrt_preprocessing.csv};

\end{groupplot}

\node at ($(group c1r1.south) + (3,-1cm)$) {\ref{commonlegend2}};
\end{tikzpicture}
\caption{Importance of preprocessing: dependency of the prediction quality $\|u-u_\NN\|$
(for the test data) on the size of the training data set and during the training}
\label{fig:performance_preprocessing}
\end{figure}

\autoref{fig:performance_preprocessing} shows the results where $y$-axis corresponds to the average values of the validation metric. As it was mentioned above, validation metric is defined as the  $\ell_2$-norm \eqref{eq:l2_norm} of difference between given and a reference solution. $x$-axis here corresponds to the number of training data.

While the topmost line depicts the validation metric of coarse (input) solutions on the test data set, the bottom-most line depicts the validation metric for the fine (target) solutions on both the test data sets. The rest of the lines depicts the validation metric for solutions obtained by using the proposed method with aforementioned preprocessing methods on train and test data sets. 

We observe that no preprocessing gives the worst results compared to other methods and the error of the resulting solution is even getting bigger with more training data.
% The second best preprocessing method turned out to be min-max scaling\comment[id=JP]{Welche Linie im Plot korrespondiert zu diesem Scaling? Wird das min-max Scaling im Text irgendwo erkl\"art.}. It performs pretty well even with small amount of the training data and the error improves with more data.

Of the preprocessing methods we tested, standardisation performed best. In this case, the error of the proposed method reaches the error of the target one slightly faster and can give the comparable error to min-max scaling while using less data for the training. Since standardization performed the best, we have used it in all our experiments. Due to the fact that standardization consists of just a few arithmetic operations performed on both input and output features, the computational overhead of it is negligible.

\subsection{Generalizations to different domains} 
In this section we discuss the application of the network to a domain different from the one, that was used for generation of the training data. In the previous sections we have used $\Omega_{\text{train}} = (0,1)^2$ for both training and testing. Now, we leave the experimental settings mostly untouched and after the training was performed on $\Omega_{\text{train}}$ we apply the network to the solutions of \eqref{eq:poisson equation} on a cross shape domain consisting of five unit squares.
In this section we consider a problems of the aforementioned form for our further tests. The equation we approximate solutions of is again the Poisson problem with the same right hand side as previously explained in \autoref{subsec:setup}.

\begin{figure}
  \begin{center}
    \includegraphics[width=.45\linewidth]{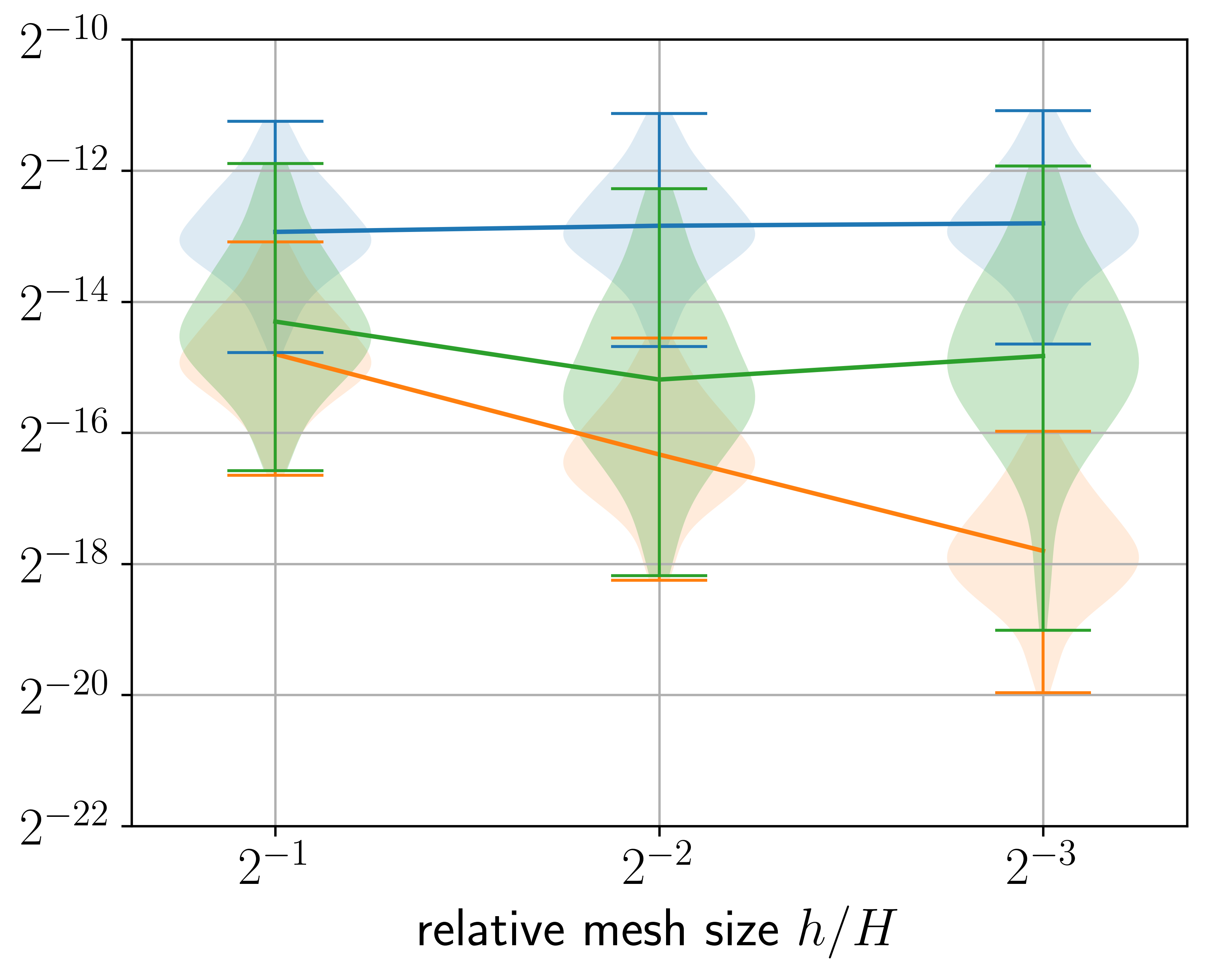}
    \caption{Distribution of the errors in the generalized domain test set for different refinement levels}
\end{center} 
\label{fig:generalization_error_distribution}
\end{figure}

\begin{figure}
  \begin{center}
    \includegraphics[width=0.21\linewidth]{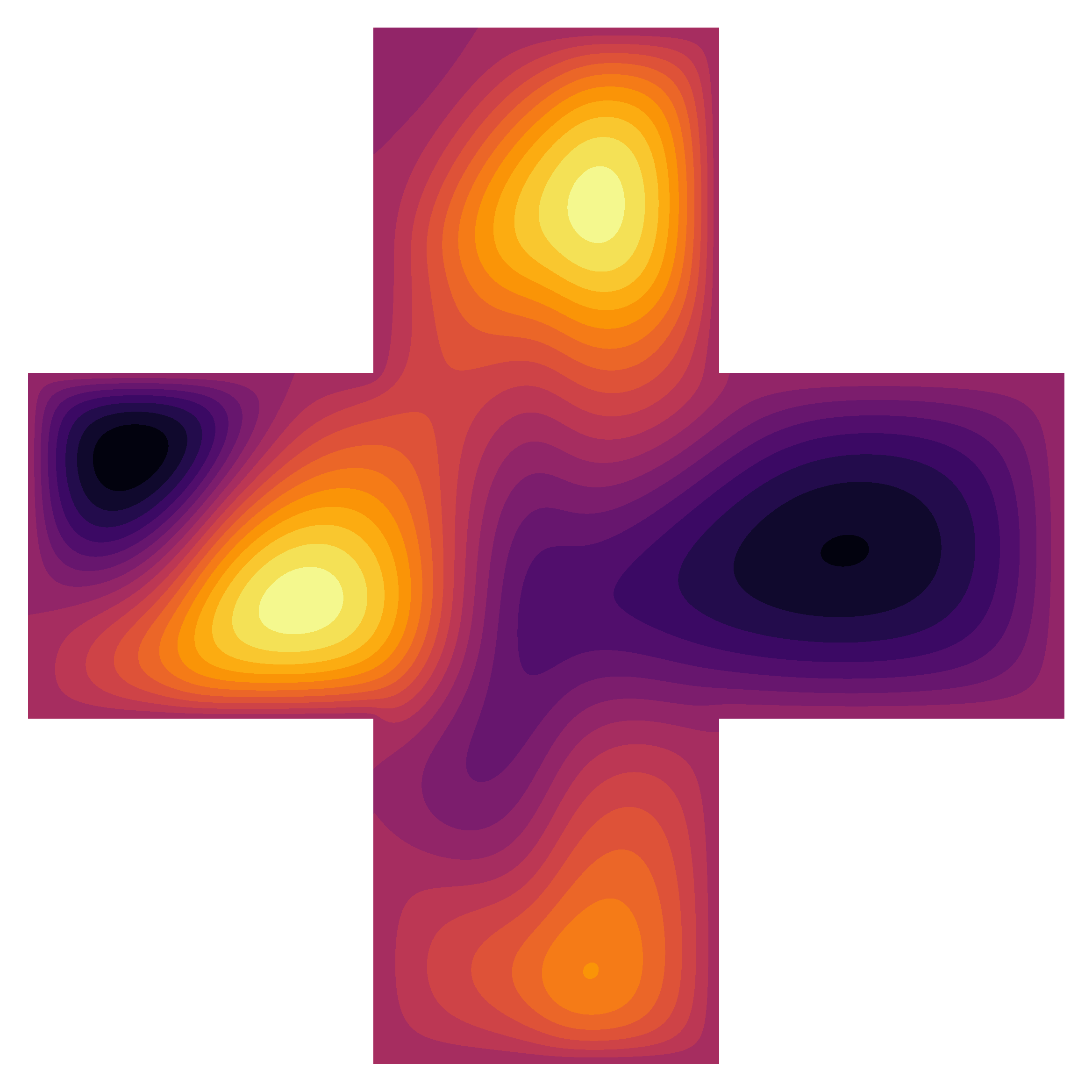}
    \includegraphics[width=0.21\linewidth]{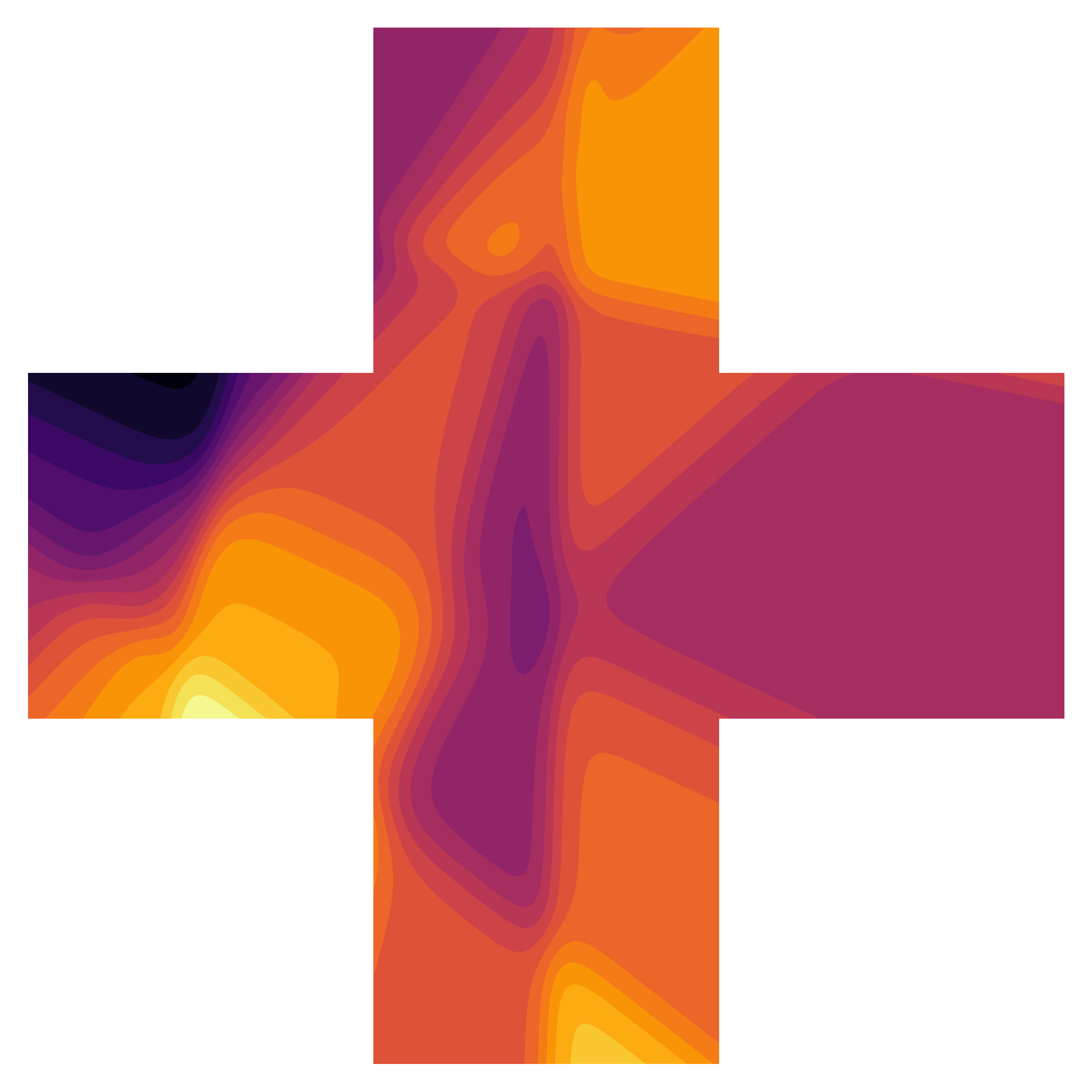}
\includegraphics[width=0.21\linewidth]{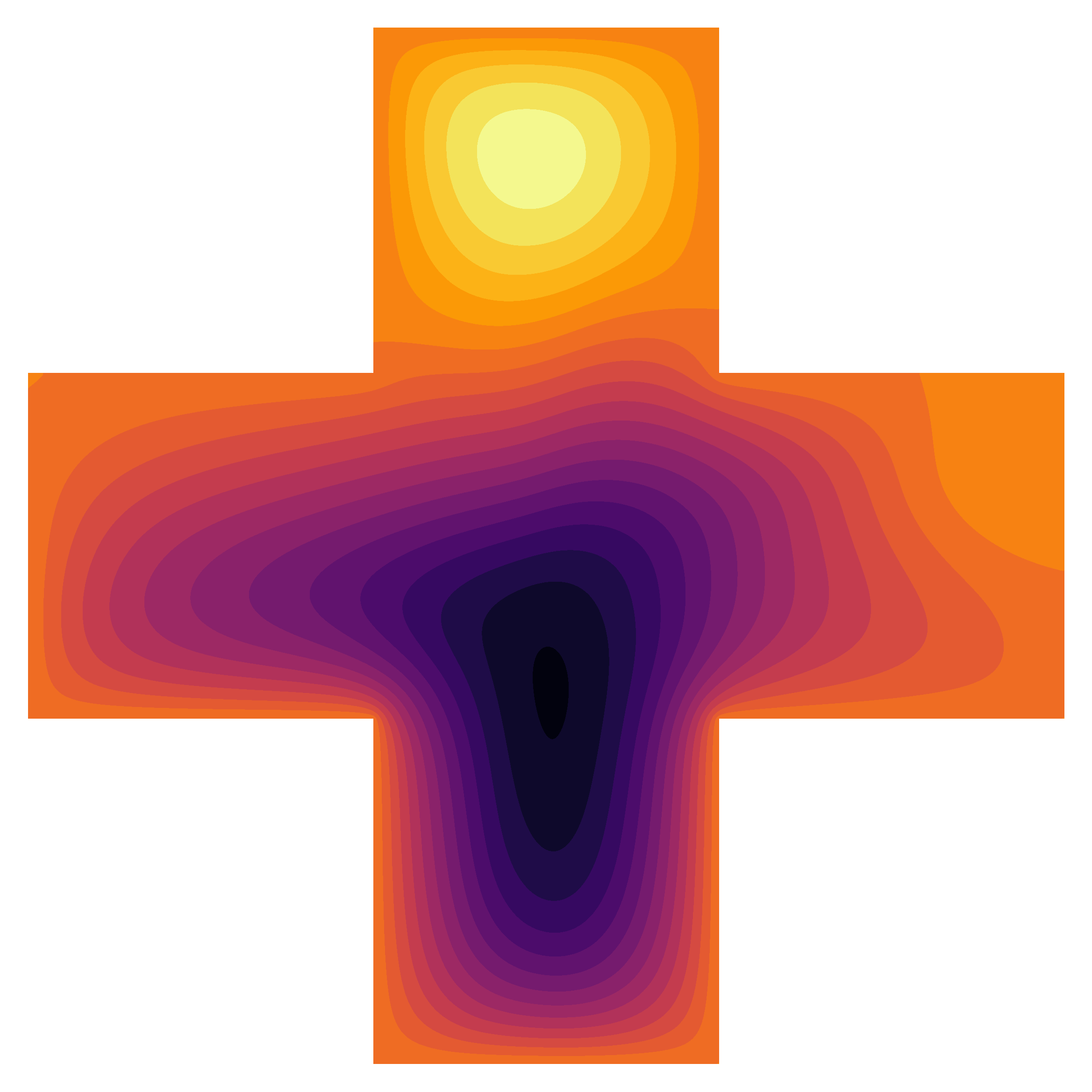}
    \includegraphics[width=0.21\linewidth]{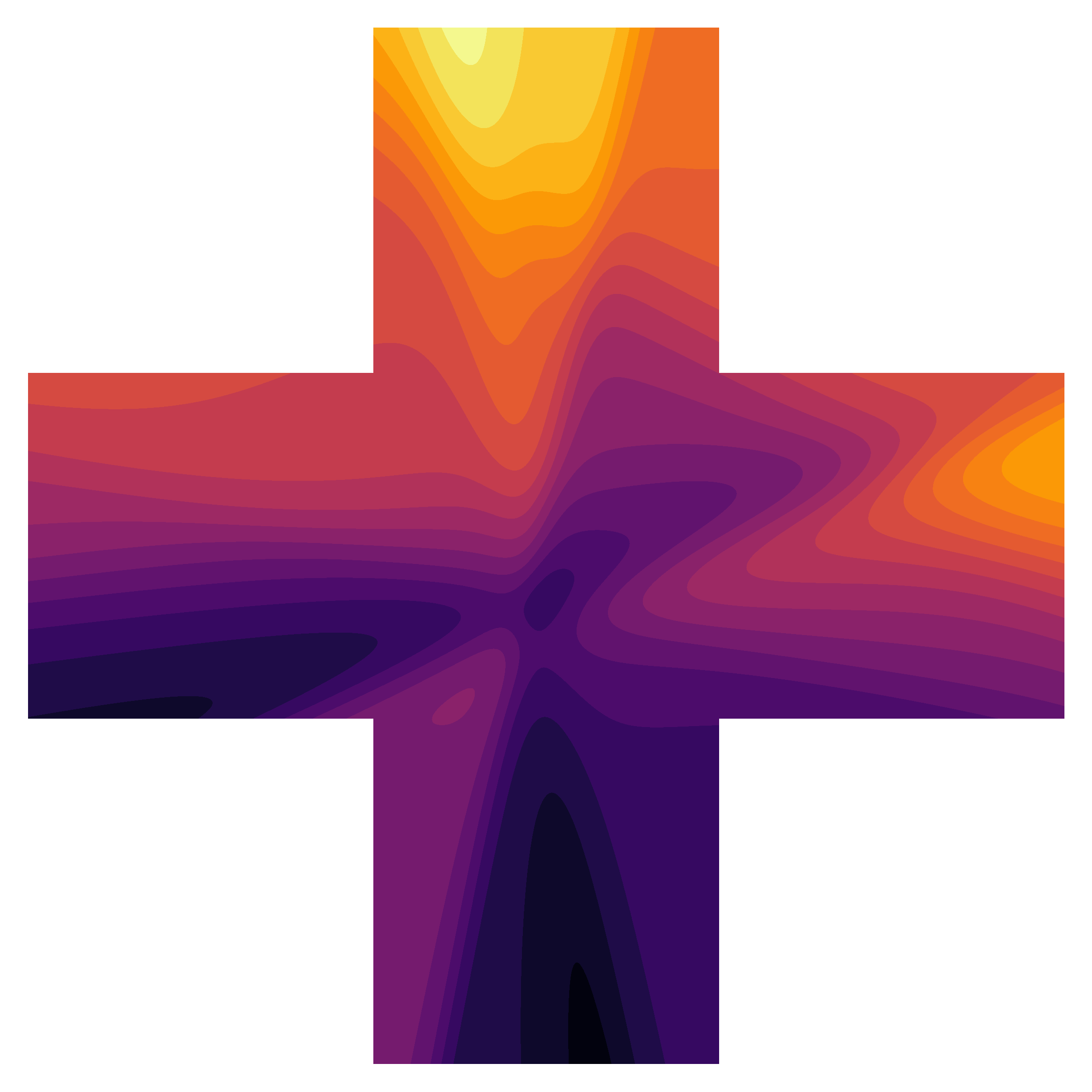}
\end{center}
\caption{Example solutions and corresponding source terms for a cross-shaped domain. Left two plots: good generalization case (solution and source term); right two plots: poor generalization case (solution and source term).}
\label{fig:generalization_solution_examples}
\end{figure}

\begin{figure}
  \begin{center}
    \includegraphics[width=0.9\linewidth]{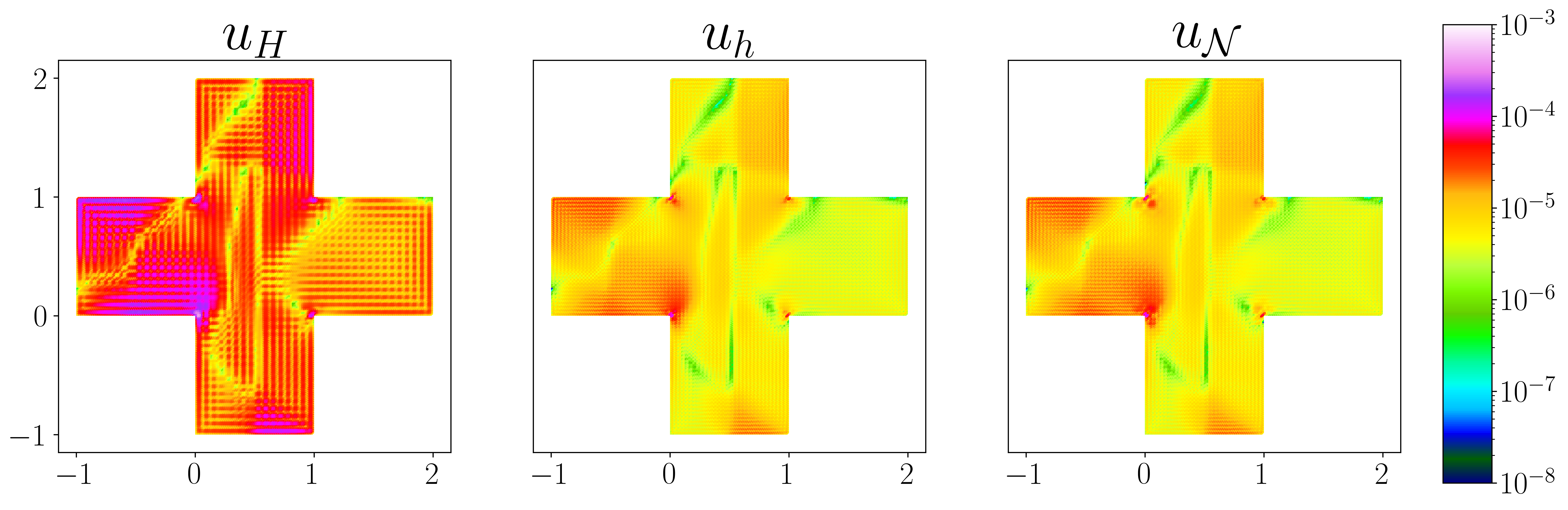}
\end{center}
    \caption{Pointwise absolute error of the coarse, fine and network solutions for a good generalization case}
\label{fig:generalization_good_errors}
\end{figure}

\begin{figure}
  \begin{center}
    \includegraphics[width=0.9\linewidth]{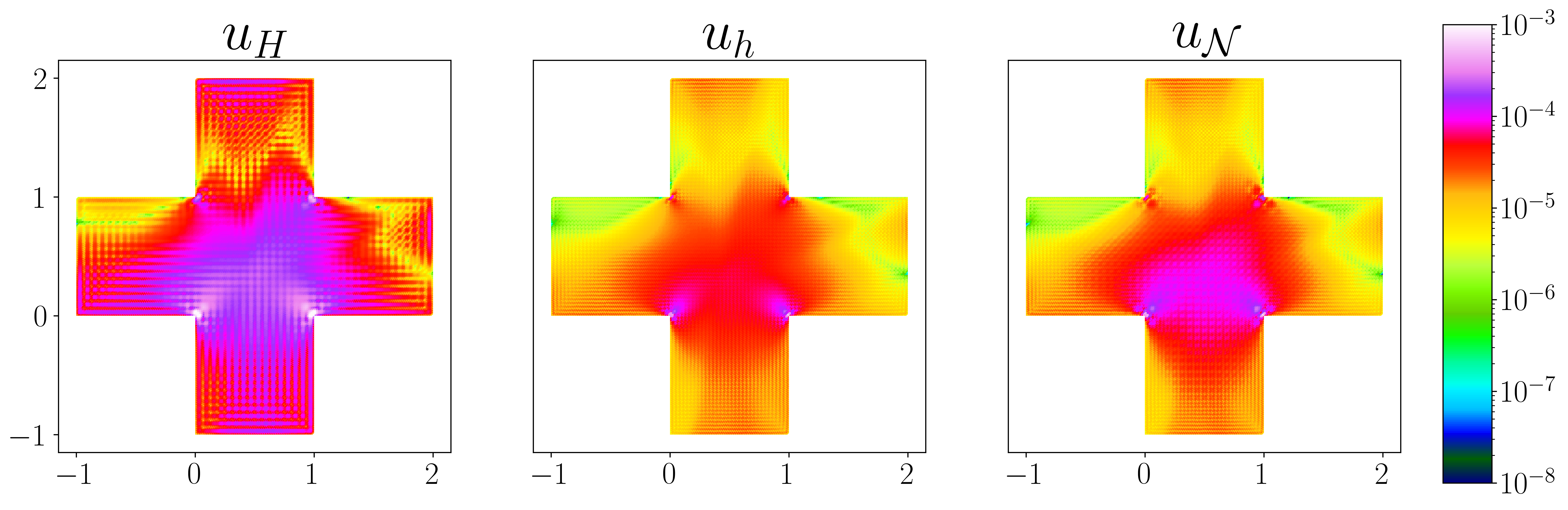}
  \end{center}
    \caption{Pointwise absolute error of the coarse, fine and network solutions for a poor generalization case}
  \label{fig:generalization_bad_errors}
\end{figure}

Table~\ref{fig:gen} shows the results for this test-case and  Figure~\ref{fig:generalization_error_distribution} shows the distribution of errors over generalization test set for different refinement levels. We observe optimal convergence for predicting one  mesh level (i.e. $h=H/2$), but the accuracy is non-optimal if two or more mesh levels are to be predicted. Figures~\ref{fig:generalization_solution_examples},  \ref{fig:generalization_good_errors} and \ref{fig:generalization_bad_errors} show the solutions and corresponding absolute pointwise error for the extended domain problem and  Table~\ref{table:generalization_errors} shows average pointwise errors for both examples.

\begin{table}
  \begin{center}
    \begin{minipage}{0.45\textwidth}
    \begin{tabular}{cccc}
      \toprule
      $h$ & $\norm{u_\text{ref} - u_H}$  & $\norm{u_\text{ref} - u_h}$  & $\norm{u_\text{ref} - u_\NN}$  \\
      \midrule
      $H/2$ & 1.46e-04 & 4.02e-05 & 6.39e-05 \\
      $H/4$ & 1.56e-04 & 1.41e-05 & 3.97e-05 \\
      $H/8$ & 1.60e-04 & 5.16e-06 & 5.17e-05 \\
      \bottomrule
    \end{tabular}
    \end{minipage}
    \begin{minipage}{0.44\textwidth}
      ~\hspace{-5.0cm}\vspace{-0.5cm}
      \begin{tikzpicture}[scale=1]
        \begin{axis}[
            width=0.9\textwidth,
            height=0.5\textwidth,
            legend style ={at={(-0.1,-0.5)},anchor=north east},
            legend columns=-1,
            xlabel={Relative mesh size $h/H$},
            ylabel={Error},
            xmode=log,
            log basis x={2},
            ymode=log,
            log basis y={2},
            mark options={solid, scale=1.6}
          ]
          \addplot[very thick, dotted, red, mark=*, mark options={solid, scale=1.2}] coordinates {(0.5,1.46e-04) (0.25,1.56e-4) (0.125,1.60e-4)};
          \addplot[very thick, dashed, blue,  mark=Mercedes star] coordinates {(0.5,4.02e-5) (0.25,1.41e-5) (0.125,5.16e-6)};
          \addplot[very thick, olive, mark=Mercedes star flipped] coordinates {(0.5,6.39e-5) (0.25,3.97e-5) (0.125,5.16e-5)};
          \legend{
	    $\|u-u_H\|$,
            $\|u-u_h\|$,
            $\|u-u_\NN\|$
          } 
        \end{axis}
      \end{tikzpicture}
    \end{minipage}
  \end{center}
  \vspace{0.5cm}
  \caption{Comparison of average errors of coarse, fine and network solutions respectively for the test problems on the generalized domain}
  \label{fig:gen}
\end{table}

% {\color{red}VLAD: can you turn this into a table-environment and add a reference to the number?}
\begin{table}
\centering
\begin{tabular}{cccc}
      \toprule
      $h$ & $\norm{u_\text{ref} - u_H}$  & $\norm{u_\text{ref} - u_h}$  & $\norm{u_\text{ref} - u_\NN}$  \\
      \midrule
      Bad example & 1.31e-04 & 4.023e-06 & 3.10e-05 \\
      Good example & 4.54e-05 & 9.76e-07 & 1.89e-06 \\
      \bottomrule
\end{tabular}
\caption{Errors for two generalization examples}
\label{table:generalization_errors}
\end{table}

% \comment[id=JP]{Figure 8 wird im Text nicht erw\"ahnt.}

\section{Conclusion}

In this paper we have introduced and analyzed a local hybrid finite element neural network method and its application to the Poisson problem. We have employed neural networks to reduce the complexity of solving Poisson equation for a family of source terms. The network operates locally, which enables it to be domain agnostic, meaning that it can be applied to problems on different domains which are not included in the training data. 

The most important component of this paper is the theoretical investigation of the proposed method. In particular we performed an a priori error analysis of the local hybrid finite element neural network method including the stability analysis of the network being used. Theoretical findings are accompanied with numerical results, which have shown that the error of the hybrid solution can be controlled by the amount of the training data and by tuning the hyperparameters of the network.

We have shown first generalization results to a slightly modified domain and found good performance. However, generalization of the method will fail if the character of the solution changes. For instance, generalization from the very regular square-domain setting to an L-shaped domain with reentrant corners will not give satisfactory results, as the singular behavior close to the corner is never seen during training of the network. Since neural networks essentially perform an interpolation of the data and processes available in the training, an appropriate enrichment of the training data will be necessary in this case.

Further steps will focus on the extension of the analysis to the time-dependent case, where the hybrid DNN-MG solver has shown very high accuracy in relevant problems~\cite{Margenberg2021,Margenberg2022,Margenberg2023}. The time-dependent case is appealing as it could enable one to achieve the accuracy of the fine temporal discretizations by calculation solutions from coarse discretizations only and receiving the temporal fluctuations from the network. This will be the focus of an upcoming work.

\section*{Data availability}
The source code for reproducing the numerical test cases can be made available on request.
% The Python scripts for reproducing the numerical test cases are published on Zenodo~\cite{zenodo}. 

\section*{Acknowledgements}

Uladzislau Kapustsin acknowledges the support of GRK 2297 MathCoRe, funded by the Deutsche Forschungsgemeinschaft, Grant Number 314838170 and the support of Federal Ministry of Research, Technology and Space as part of the research grant 03VP13242.
Utku Kaya acknowledges the support of GRK 2297 MathCoRe, funded by the Deutsche Forschungsgemeinschaft, Grant Number 314838170.
Thomas Richter acknowledges the support of GRK 2297 MathCoRe, funded by the Deutsche Forschungsgemeinschaft, Grant Number 314838170, 
the support of Federal Ministry of Research, Technology and Space as part of the research grant 03VP13242 and
the support of the Deutsche Forschungsgemeinschaft, Grant Number 537063406.

\section*{Conflict of interest}
The authors declare that they have no conflict of interest.

\bibliography{references}

\end{document}